\def\date{5.9.06}   

\documentclass[12pt]{article}

\usepackage{amssymb}
\usepackage{amsmath}
\input{liemacs.sty} 

\newcommand{\epi}{\mathop{{\rm epi}}\nolimits}
\newcommand{\Rep}{\mathop{{\rm Rep}}\nolimits}

\begin{document} 



\title{A complex semigroup approach to group algebras of infinite dimensional 
Lie groups}

\author{Karl-Hermann Neeb}

\maketitle

\centerline{\sl To K. H. Hofmann on the occasion of his 75th birthday} 

\begin{abstract} A host algebra of a topological group $G$ is a $C^*$-algebra 
whose representations are in one-to-one correspondence with certain 
continuous unitary representations of $G$. In this paper we present an approach 
to host algebras for infinite dimensional Lie groups which is based on 
complex involutive semigroups. Any locally bounded absolute value $\alpha$ 
on such a semigroup $S$ leads in a natural way to a 
$C^*$-algebra $C^*(S,\alpha)$, and we describe a setting which permits us 
to conclude that this $C^*$-algebra is a host algebra for a Lie group $G$. 
We further explain how to attach to any such host algebra an invariant 
weak-$*$-closed convex set in the dual of the Lie algebra of $G$ enjoying certain 
nice convex geometric properties. If $G$ is the additive group of a locally 
convex space, we describe all host algebras arising this way. The general
non-commutative case is left for the future. 
\end{abstract} 

\section*{Introduction} 

If $G$ is a locally compact group, then Haar measure 
on $G$ leads to the convolution algebra $L^1(G)$, and we obtain 
a $C^*$-algebra $C^*(G)$ as the enveloping $C^*$-algebra of $L^1(G)$. 
This $C^*$-algebra has the universal property that each 
(continuous) unitary representation $(\pi, {\cal H})$ of $G$ on 
some Hilbert space ${\cal H}$ defines a unique non-degenerate 
representation of $C^*(G)$ on ${\cal H}$ and, conversely, each 
non-degenerate representation of $C^*(G)$ arises from a unique 
unitary representation of $G$. This correspondence is a central tool 
in the harmonic analysis on $G$ because the well-developed theory of 
$C^*$-algebras provides a powerful machinery to study the set of all 
irreducible representations of $G$, to endow it with a natural 
topology and to understand how to decompose representations into 
irreducibles or factor representations. 

For infinite dimensional Lie groups, i.e., Lie groups modeled on 
infinite dimensional locally convex spaces, there is no natural 
analog of the convolution algebra $L^1(G)$, so that we cannot hope 
to find a $C^*$-algebra whose representations are in one-to-one 
correspondence to {\sl all} unitary representations of $G$. However, 
in \cite{Gr05} H.~Grundling introduces the notion of a {\it host algebra}  
of a topological group $G$. This is a pair $({\cal A}, \eta)$, consisting 
of a $C^*$-algebra ${\cal A}$ and a morphism 
$\eta \: G \to U(M({\cal A}))$ of $G$ into the unitary group 
of its multiplier algebra $M({\cal A})$ with the following property:  
For each non-degenerate representation $\pi$ of ${\cal A}$ 
and its canonical extension $\tilde\pi$ to $M({\cal A})$, 
the unitary representation $\tilde\pi \circ \eta$ of $G$ is continuous 
and determines $\pi$ uniquely. In this sense, ${\cal A}$ is hosting a 
certain class of representations of $G$. A host algebra $\cA$ is called 
{\it full} if it is hosting all continuous unitary representations of $G$.  
Now it is natural to ask to which extent infinite dimensional 
Lie groups, or other non-locally compact groups, possess host algebras. 
One cannot expect the existence 
of a full host algebra because, f.i., the topological dual $E'$ of 
an infinite dimensional locally convex space $E$ carries no natural 
locally compact topology. Therefore one is looking for host algebras that 
accommodate certain classes of continuous unitary representations. 

In the present paper we discuss a construction of host algebras 
based on holomorphic extensions of unitary representations of a 
Lie group $G$ to certain complex semigroups $S$. Some of the basic 
ideas of our constructions appear already in \cite{Ne95}, 
where one finds the construction of the 
enveloping $C^*$-algebras $C^*(S,\alpha)$ of a complex involutive 
semigroup $S$, endowed with a locally bounded absolute value $\alpha$, 
and also in \cite{Ne98}, where this is applied 
to the special case where $S$ is a complex Banach--Lie group. Here we 
address the situation  where $S$ may be an infinite dimensional 
semigroup which is not a group. 

The structure of the paper is as 
follows. In Section 1 we first recall the concept 
of a complex involutive semigroup $S$ and associate to any locally 
bounded absolute value $\alpha$ on $S$ a $C^*$-algebra $C^*(S,\alpha)$ 
with a holomorphic morphism $\eta_\alpha \: S \to C^*(S,\alpha)$ 
having a suitable universal property. Since our goal is to construct 
host algebras for infinite dimensional Lie groups, 
we build in Section 2 a bridge between complex involutive semigroups 
and Lie groups by defining the notion of a host semigroup of a Lie group. 
Roughly speaking, this a complex involutive semigroup $S$ on which the Lie group 
$G$ acts smoothly by unitary multipliers and for which there exists an 
open convex cone $W$ in the Lie algebra $\L(G)$ of $G$, 
invariant under the adjoint action, 
for which all $\R$-actions on $S$ defined by the one-parameter semigroups 
$\gamma_x(t) = \exp_G(tx)$, $x \in W$, extend to ``holomorphic'' 
one-parameter semigroups $\C_+ = \R + i [0,\infty[ \to M(S)$ mapping 
the open upper halfplane $\C_+^0$ holomorphically into $S \subeq M(S)$.
The main result of Section 2 is that for 
each locally bounded absolute value $\alpha$ on a host semigroup $S$, the 
$C^*$-algebra $C^*(S,\alpha)$ is a host algebra of $G$. 

This leaves us with the problem to understand the classes of representations 
of $G$ hosted by such $C^*$-algebras. To clarify this point, we consider 
in Section 3 multiplier actions $\eta \: G \to U(M(\cA))$ 
of a Lie group $G$ on a $C^*$-algebra 
$\cA$ and study to which extent the action of certain one-parameter 
semigroups of $G$ extends holomorphically to the upper halfplane. 
This leads  to the momentum map 
$$ \Psi_\eta \: S(\cA)^\infty \to \L(G)', \quad 
\phi \mapsto \frac{1}{i} d(\tilde\phi \circ \eta)(\1). $$
Here $S(\cA)^\infty$ denotes the set of all states $\phi$ of $\cA$ 
for which the canonical extension $\tilde\phi$ to $M(\cA)$ yields a 
smooth function $\tilde\phi \circ \eta \: G \to \C$. 
The  weak-$*$-closed hull $I_\eta$ of the image of $\Psi_\eta$ 
is a convex set invariant under 
the coadjoint action, called the {\it momentum set of $(\cA,\eta)$}. 
A crucial observation is that, if the multiplier 
action comes from a host algebra $C^*(S,\alpha)$, where $S$ is a host semigroup, 
the convex cone 
$$ B(I_\eta) := \{ x \in \L(G) \: \inf \la I_\eta, x \ra > - \infty \} $$
has non-empty interior and the support function 
$$ s \: B(I_\eta)^0 \to \R, \quad x\mapsto - \inf \la I_\eta, x \ra $$
is locally bounded. This observation suggests that to find host algebras for 
$G$, one should start with an $\Ad^*(G)$-invariant weak-$*$-closed convex subset 
$C \subeq \L(G)'$ for which the corresponding {\it support function} 
$s_C \: B(C)^0 \to \R$ 
is locally bounded. As the function $s_C$ is convex, we take in 
Section 4 a closer look at convex functions on open convex domains in locally 
convex spaces. In particular, we show that whenever 
$\L(G)$ is barrelled, the existence of interior points in the cone 
$B(C)$ automatically implies that $s_C$ is locally bounded and even continuous. 

In Section 5 we then show how this circle of ideas can be completed in 
the abelian case. Here the Lie group $G$ is a locally convex space 
$V$ and the semigroup is a tube $S = V + i W$, where 
$W \subeq V$ is an open convex cone. In this case it suffices to 
consider absolute values of the form \break 
$\alpha(a + ib) = e^{-\inf \la C, b \ra}$, where 
$C \subeq V'$ is a weak-$*$-closed convex subset. 
Now $\alpha$ is locally bounded if and only if the support function 
$s_C(x) = - \inf \la C, x \ra$ is locally bounded on $W$. 
If this is the case, then the results of Section 4 imply that 
$C$ is locally compact and $C^*(S,\alpha) \cong C_0(C)$ is a host algebra 
of $V$ hosting precisely all unitary representations of $V$ arising from 
spectral measures on the locally compact subset $C\subeq V'$. 

Section 6 contains a brief discussion of the finite dimensional case, 
which is developed in detail in \cite{Ne99}. Here we give a short 
and direct proof of the fact that any host algebra of $G$ coming from a 
host semigroup is a quotient of the group $C^*$-algebra $C^*(G)$. 

The next steps of this project aim at a better understanding of the classes of 
representations of a Lie group $G$ hosted by $C^*$-algebras of the form 
$C^*(S,\alpha)$. The first major problem one has to solve here 
is to find a suitable complex involutive semigroup $S$ 
whenever the invariant convex set $C \subeq \L(G)'$ is given. For 
finite dimensional 
groups this has been carried out in \cite{Ne99}, but for infinite dimensional 
groups many key tools are still missing. Furthermore, once the semigroup $S$ is 
constructed, one has to find the class of unitary representations of $G$ 
extending to holomorphic representations of $S$. We leave all that 
to the future. 

\subsection*{Preliminaries} 

For the sake of easier reference, we collect some of the basic definitions 
concerning infinite dimensional manifolds and Lie groups. 

Let $X$ and $Y$ be locally convex topological vector spaces, $U
\subeq X$ open and $f \: U \to Y$ a map. Then the {\it derivative
  of $f$ at $x$ in the direction of $h$} is defined as 
$$ df(x)(h) := \lim_{t \to 0} \frac{1}{t} \big( f(x + t h) - f(x)\big)
$$
whenever the limit exists. The function $f$ is called {\it differentiable at
  $x$} if $df(x)(h)$ exists for all $h \in X$. It is called {\it
  continuously differentiable or $C^1$} if it is continuous and differentiable at all
points of $U$ and 
$$ df \: U \times X \to Y, \quad (x,h) \mapsto df(x)(h) $$
is a continuous map. It is called a $C^n$-map if it is $C^1$ and $df$ is a
$C^{n-1}$-map, and $C^\infty$ (or {\it smooth}) if it is $C^n$ for all $n \in \N$. 
This is the notion of differentiability used in \cite{Mil84}, \cite{Ha82} and
\cite{Gl02}, where the latter reference deals with the modifications
needed for incomplete spaces. If $X$ and $Y$ are complex, 
$f$ is called {\it holomorphic} if it is smooth 
and its differentials $df(x)$ are 
complex linear. If $Y$ is Mackey complete, it suffices that $f$ is $C^1$. 

Since we have a chain rule for $C^1$-maps between locally convex 
spaces, we can define smooth manifolds as in
the finite dimensional case. A chart $(\phi,U)$ with respect to a given 
manifold structure on $M$ is an open set $U\subset M$ together with a 
homeomorphism $\phi$ onto an open set of the model space. 

A {\it Lie group} $G$ is a group equipped with a 
smooth manifold structure modeled on a locally convex space 
for which the group multiplication and the 
inversion are smooth maps. We write $\1 \in G$ for the identity element and 
$\lambda_g(x) = gx$, resp., $\rho_g(x) = xg$ for the left, resp.,
right multiplication on~$G$. Then each $x \in T_\1(G)$ corresponds to
a unique left invariant vector field $x_l$ with 
$x_l(\1) = x$. 
The space of left invariant vector fields is closed under the Lie
bracket, hence inherits a Lie algebra structure. We thus obtain 
on the tangent space $T_\1(G)$ a continuous Lie bracket which
is uniquely determined by $[x,y]_l = [x_l, y_l]$ for $x,y \in T_\1(G)$. 
We write $\L(G) := (T_\1(G), [\cdot,\cdot])$ for the so-obtained 
topological Lie algebra. Then $\L$ defines a functor from the category 
of locally convex Lie groups to the category of locally convex topological 
Lie algebras. The adjoint action of $G$ on $\L(G)$ is defined by 
$\Ad(g) := \L(c_g)$, where $c_g(x) = gxg^{-1}$. This action is smooth and 
each $\Ad(g)$ is a topological isomorphism of $\L(G)$. The coadjoint action 
on the topological dual space $\L(G)'$ is defined by $\Ad^*(g).f := f \circ 
\Ad(g)^{-1}$ and all these maps are continuous with respect to the 
weak-$*$-topology on $\L(G)'$, but in general the coadjoint 
action of $G$ is {\sl not} continuous with respect to this topology. 

\section{$C^*$-algebras associated to complex semigroups} 

In this section we associate to each complex involutive semigroup 
$S$, endowed with an absolute value $\alpha$, a $C^*$-algebra $C^*(S,\alpha)$. 
As we shall see later on, one can use these semigroup algebras to construct 
host algebras for Lie groups, and this is our main main purpose 
for their construction.

\begin{definition}
  \label{def:1.1} 
(a) An {\it involutive complex semigroup} is 
a complex manifold $S$ modeled on a locally convex space which is endowed 
with a holomorphic semigroup multiplication and an antiholomorphic 
antiautomorphism denoted $s \mapsto s^*$. 

\par (b) A function $\alpha \: S \to \R^+$ 
is called an {\it absolute value} if 
$$ \alpha(s) = \alpha(s^*) \quad \mbox{ and } \quad 
\alpha(st) \leq \alpha(s) \alpha(t) $$
for all $s,t \in S$. 

\par (c) A {\it holomorphic representation} $(\pi, {\cal H})$ 
of a complex involutive semigroup $S$ on the Hilbert 
space ${\cal H}$ is a morphism $\pi \: S \to B({\cal H})$ 
of involutive semigroups which is 
holomorphic if $B({\cal H})$ is endowed with its natural complex Banach
space structure defined by the operator norm. 
If $\alpha$ is an absolute value on $S$, then the 
representation $\pi$ is said to be {\it $\alpha$-bounded} if 
$\|\pi(s)\| \leq \alpha(s)$ holds for each $s \in S$. 
A representation is called {\it non-degenerate} if 
$\pi(S).v = \{0\}$ implies $v= 0$. 
\end{definition}

\begin{examples} (1) If $H$ is a complex Lie group and $s \mapsto s^*$ an 
antiholomorphic antiautomorphism $H$, then $H$ is a complex involutive (semi)group. 
Any open $*$-subsemigroup of $H$ is a complex involutive semigroup. 

(2) If $V$ is a locally convex space and $W \subeq V$ an open convex cone, then 
$S := V + i W \subeq V_\C$ is an involutive subsemigroup with respect to the 
involution $(x + iy)^* := -x + iy$.  

(3) If $\cA$ is a $C^*$-algebra, then its multiplicative semigroup $(\cA,\cdot)$ 
is a complex involutive semigroup and $\alpha(a) := \|a\|$ is an absolute value 
on $\cA$.   
An important example is the $C^*$-algebra $B({\cal H})$ of bounded operators on 
the Hilbert space ${\cal H}$. 

(4) Let $\cA$ be a unital $C^*$-algebra and $\tau = \tau^* \in \cA$ an 
involution, i.e., $\tau^2 = \1$. 
For $a,b \in \cA$ we write $a < b$ if there exists an invertible element $c \in \cA$ with $b - a = c^*c$. Then 
$$ S:= \{ s \in A \: s^*\tau s < \tau \} $$
is an open subsemigroup of $\cA$ with respect to multiplication. 
To see that it is non-empty, we observe that we may write 
$\tau = \1 - 2 p = (\1 - p) - p$ 
for a projection $p = p^* = p^2 \in \cA$. For $\lambda \in \C^\times$ and 
$s := \lambda (\1-p) + \lambda^{-1}p$ we then have 
$$ s^* \tau s 
= |\lambda|^2 (\1 - p) - |\lambda^{-1}|^2 p < \tau = (\1-p) - p $$
if and only if $|\lambda| < 1$. The boundary of $S$ contains 
the real Banach--Lie group 
$U(\cA,\tau) := \{ g \in \cA^\times \: g^*\tau g = \tau \}.$
\end{examples}

\begin{definition} Let $S$ be a complex involutive semigroup 
and $\alpha$ a locally bounded absolute value on $S$. We associate 
to the pair $(S,\alpha)$ a $C^*$-algebra $C^*(S,\alpha)$ as follows. 

First, we endow the semigroup algebra $\C[S]$, whose 
elements we write as finitely supported functions $f \: S \to \C$, 
with the submultiplicative seminorm 
$\|f\|_\alpha := \sum_{s\in S} |f(s)|\alpha(s)$ 
and the involution $f^*(s) := \oline{f(s)}$. 
Let $\ell^1(S,\alpha)$ be the complex involutive Banach algebra obtained 
by completion of this seminormed $*$-algebra. 
We define $\eta_\alpha^1(s) \in \ell^1(S,\alpha)$ as the image 
of the function $\delta_s(t) := \delta_{s,t}$ in $\ell^1(S,\alpha)$ and 
note that $\|\eta_\alpha(s)\| = \alpha(s)$. 

If $\cA$ is a $C^*$-algebra, then each homomorphism 
$\beta \: S \to (\cA,\cdot)$ of involutive semigroups, which is $\alpha$-bounded 
in the sense that $\|\beta(s)\| \leq \alpha(s)$ holds for each $s \in S$,  
defines a unique contractive morphism 
$$\hat\beta \: \ell^1(S,\alpha) \to \cA, \quad \hat\beta(f) := \sum_{s \in S} 
f(s)\beta(s) $$ 
of Banach-$*$-algebras satisfying $\hat\beta \circ \eta_\alpha^1 = \beta$. 
Let $I \trile \ell^1(S,\alpha)$ denote the 
intersection of the kernels of all such homomorphism $\hat\beta$ for which 
$\beta$ is a holomorphic map. On the quotient algebra 
$\ell^1(S,\alpha)/I$, we obtain a $C^*$-norm by 
$$ \|[f]\| := \sup_{\beta\ {\rm holomorphic}} \|\hat\beta(f)\| \leq \|f\|_\alpha. $$
We now define $C^*(S,\alpha)$ as the completion of 
$\ell^1(S,\alpha)/I$ with respect to this norm. It follows immediately 
from the construction that we thus obtain a $C^*$-algebra. 
\end{definition}

Before we turn to the universal property of $C^*(S,\alpha)$, 
we recall the following criterion for holomorphy (\cite{Ne99}, Cor.~A.III.3): 

\begin{lemma} \label{lem:holcrit} 
Let $M$ be a complex manifold, $V$ a Banach space
and $N \subeq V'$ a subset which is norm-determining, i.e., 
$\|v\| = \sup \{ |\lambda(v)| 
\: \lambda \in N, \|\lambda\| \leq~1\}$ for all $v \in V$. Then a
locally bounded function $f \: M \to V$ is holomorphic if and only if 
for each $\lambda \in N$ the function $\lambda \circ f$ is holomorphic. 
\end{lemma}

The following theorem could also be derived from 
Theorem~IV.2.7 in \cite{Ne99}, but the construction 
we give here is much more direct. 

\begin{theorem}
  \label{thm:1.7} 
The $C^*$-algebra 
$C^*(S,\alpha)$ has the following properties: 
\begin{description}
\item[\rm(i)] There exists a holomorphic morphism 
$\eta_\alpha\: S \to C^*(S,\alpha)$ of in\-vo\-lu\-tive semigroups with total
range, i.e., $\eta_\alpha(S)$ generates a dense subalgebra. 
\item[\rm(ii)] For each $\alpha$-bounded holomorphic morphism of involutive semigroups \break $\pi \: S \to {\cal A}$ to the multiplicative semigroup of 
a $C^*$-algebra ${\cal A}$, there exists a unique 
morphism of $C^*$-algebras $\tilde\pi \: C^*(S,\alpha) \to \cA$ 
with $\tilde\pi \circ \eta_\alpha = \pi$. 
\end{description}
\end{theorem}

\begin{proof} (i) We define $\eta_\alpha(s) \in C^*(S,\alpha)$ as the image 
of the element $\eta_\alpha^1(s) \in \ell^1(S,\alpha)$. 
Then $\|\eta_\alpha(s)\| \leq \|\eta_\alpha^1(s)\| = \alpha(s)$ implies that 
$\eta_\alpha$ is a locally bounded morphism of involutive semigroups. 

To see that $\eta_\alpha$ is holomorphic, we first note that the subspace 
$N$ of continuous linear functionals on $C^*(S,\alpha)$ spanned by 
the functionals of the form $\phi \circ \hat\beta$ on 
$\ell^1(S,\alpha)$, where $\beta \: S \to \cA$ 
is a holomorphic morphism of involutive semigroups into a $C^*$-algebra $\cA$ 
and $\phi \in {\cal A}'$,  
separates the points of $C^*(S,\alpha)$ and determines the norm 
(by definition of the norm on $C^*(S,\alpha)$). 
For each functional $\psi = \phi \circ \hat\beta$ as above 
the map 
$$ \psi \circ \eta_\alpha = \phi \circ \beta \: S \to \C $$
is holomorphic. 
Therefore Lemma~\ref{lem:holcrit} implies that $\eta_\alpha$ is holomorphic. 

That $\eta_\alpha(S)$ spans a dense subspace of $C^*(S,\alpha)$ follows from the 
construction because the image of $S$ spans a dense subspace of $\ell^1(S,\alpha)$, 
hence also in the quotient by the ideal $I$. 

(ii) Let $\tilde\pi \: \ell^1(S,\alpha) \to \cA$ denote the canonical 
extension of $\pi$ which is a contractive morphism of involutive Banach 
algebras and note that $\ker \tilde\pi \supeq I$, so that 
$\tilde\pi$ factors through a morphism $\ell^1(S,\alpha)/I \to \cA$ of 
involutive Banach algebras which, by definition, extends to the completion 
$C^*(S,\alpha)$. 
\end{proof}

\begin{remark} \label{rem:adjoint} (a) The 
preceding theorem entails that for each $C^*$-algebra $\cA$, we have 
$\Hom(C^*(S,\alpha),\cA) \cong \Hom_{\rm hol}((S,\alpha),(\cA,\|\cdot\|)$, where 
the right hand side denote the holomorphic contractive 
morphisms of complex involutive semigroups 
with absolute value. This means that $C^*(S,\alpha)$ defines an adjoint of the 
forgetful functor from the category of $C^*$-algebras to the category of 
complex involutive semigroups with absolute value, assigning to a $C^*$-algebra 
$\cA$ the semigroup $(\cA,\cdot, \|\cdot\|)$. 
It follows in particular that the universal property in 
Theorem~\ref{thm:1.7}(ii) determines $C^*(S,\alpha)$ up to isomorphism. 

(b) If $\beta$ is an absolute value on $S$ satisfying 
$\|\eta_\alpha\| \leq \beta \leq \alpha$, then the natural map 
$\phi \: C^*(S,\beta) \to C^*(S,\alpha)$ is an isomorphism with 
$\phi \circ \eta_\beta = \eta_\alpha$. 
\end{remark}

\begin{examples} \label{ex:cstarsem} (a) We 
take a closer look at the case where $S$ is commutative. 
Let $\hat S:= \Hom(S,(\C,\cdot)) \setminus\{0\}$ denote the set of 
non-zero holomorphic characters of $S$, i.e., the one dimensional 
(=irreducible) non-degenerate representations. A holomorphic 
character 
$\chi$ extends to a character of the $C^*$-algebra $C^*(S,\alpha)$ if and only 
if it is $\alpha$-bounded. Hence the set $\hat S_\alpha$ of $\alpha$-bounded 
non-zero holomorphic characters form the spectrum of the commutative 
$C^*$-algebra $C^*(S,\alpha)$. We conclude that 
$C^*(S,\alpha) \cong C_0(\hat S_\alpha)$, where 
$\hat S_\alpha \subeq C^*(S,\alpha)'$ carries the weak-$*$-topology. 
Moreover, the set $\hat S_\alpha \cup \{0\}$ is weak-$*$-compact in 
$C^*(S,\alpha)'$, and the canonical map 
$\eta_\alpha^* \: 
C^*(S,\alpha)' \to \C^S, \phi \mapsto \phi \circ \eta$ 
is continuous with respect to the weak-$*$-topology 
on the left and the product topology on the right. 
This shows that $\hat S_\alpha \cup \{0\}$ is compact in $\C^S$ with respect to 
the product topology which, therefore, coincides with 
the weak-$*$-topology defined by $C^*(S,\alpha)$. 
We conclude that $\hat S_\alpha$ is locally compact with respect to the 
product topology and that $C^*(S,\alpha) 
\cong C_0(\hat S_\alpha)$. We now have 
$$ \|\eta_\alpha(s)\| = \sup \{ |\chi(s)| \: \chi \in \hat S_\alpha\}, $$
and this absolute value defines the same $C^*$-algebra 
by Remark~\ref{rem:adjoint}(b). 

(b) We specialize to the particular case where 
$S = V + i W \subeq V_\C$ holds for a real 
locally convex space $V$ and an open 
convex cone $W \subeq V$. This is an open complex subsemigroup of 
the complex vector space $V_\C$. Any non-zero holomorphic character 
$\chi \: S \to \C$ maps into $\C^\times$ and induces a unique continuous 
character $V \to \T$, hence is of the form 
$\chi = e^{if}$ for some continuous linear functional $f \in V'$ 
(which we also extend to a complex linear functional on $V_\C$). 

Now let $\alpha$ be a locally bounded 
absolute value on $S$ and 
$$ C_\alpha := \{ f \in V' \: e^{if} \in \hat S_\alpha \} $$ 
the set of linear functionals defining 
$\alpha$-bounded characters of $S$. In view of (b), we may 
w.l.o.g.\ assume that 
$$\alpha(x+iy) = \|\eta_\alpha(x+iy)\| 
= \sup \{ e^{-f(y)} \: f \in C_\alpha\}
= e^{-\inf \la C_\alpha, y \ra} $$
without changing $C^*(S,\alpha)$ or $C_\alpha$. 
A holomorphic character $e^{if}$, $f \in V'$, 
is $\alpha$-bounded, i.e., contained in $\hat S_\alpha$, if and only if
$e^{-f(y)} \leq e^{-\inf \la C_\alpha, y \ra}$
for each $y \in W$, which is equivalent to 
\begin{equation}
  \label{eq:sstar}
f(y) \geq \inf \la C_\alpha, y \ra \quad \mbox { for } \quad y \in W. 
\end{equation}
This implies that $C_\alpha$ is a weak-$*$-closed convex subset 
of $V'$, and (a) further shows that $\hat S_\alpha$ is locally compact 
as a subset of $\C^S$. We continue the discussion of these examples in 
Section 5 below. 

(c) Let $S := \C_+^0 
= \R + i ]0,\infty[ \subeq \C$ be the open upper half plane 
and $\alpha$ as in (b). Then $V = \R$, $W = ]0,\infty[$ and 
$C_\alpha \subeq \R$ is a closed convex subset bounded from below. 
Let $m := \inf C_\alpha$. Then 
$$ \|\eta_\alpha(x + iy)\| = e^{-(\inf C_\alpha) \cdot y} = e^{-my}, $$
and (\ref{eq:sstar}) implies that $C_\alpha = [m,\infty[$. 
Therefore $C^*(\C_+^0, \alpha) \cong C_0([m,\infty[)$. 
\end{examples}

\subsection*{Some holomorphic representation theory}

\begin{lemma}
  \label{lem:1.2} 
Let ${\cal H}$ be a Hilbert space and 
$\pi \: S \to B({\cal H})$ a morphism of involutive semigroups. 
Then $\pi$ is a holomorphic representation 
if and only if it satisfies the following conditions: 
\begin{description}
\item[\rm(1)]$\pi$ is locally bounded, i.e., for every $s \in S$ there exists a 
neighborhood $U$ such that $\pi(U)$ is a bounded subset of the Banach space 
$B({\cal H})$. 
\item[\rm(2)] There exists a dense subspace $E \subeq {\cal H}$ such 
that the functions $\pi^v(s) := \la \pi(s).v,v\ra$ 
are holomorphic for all $v \in E$.
\end{description}
\end{lemma}

\begin{proof}
The necessity of conditions (1) and (2) is obvious, and the
converse follows from \cite{Ne99}, Cor.~A.III.5, which also holds for 
general locally convex manifolds since \cite{He89}, Prop.~2.4.9(a) applies 
to functions on locally convex spaces that are not necessarily complete. 
\end{proof}

\begin{proposition}
  \label{prop:1.3} Let $S$ be an involutive complex semigroup and 
$\alpha$ a locally bounded absolute value on $S$. 

\begin{enumerate}
\item[\rm(a)] If $(\pi_j, {\cal H}_j)_{j \in J}$ is a set of $\alpha$-bounded holomorphic 
representations of $S$, then the operators induced by 
$s \in S$ on $\bigoplus_{j \in J} {\cal H}_j$ are 
bounded, and we thus obtain an $\alpha$-bounded 
holomorphic representation of $S$ on 
the direct Hilbert space sum $\hat{\bigoplus}_{j \in J} {\cal H}.$ 

\item[\rm(b)] Every non-degenerate $\alpha$-bounded holomorphic representation 
is a direct sum of cyclic $\alpha$-bounded holomorphic representations.  

\item[\rm(c)] Every $\alpha$-bounded cyclic holomorphic representation of $S$ is
equivalent to a representation $(\pi_\phi, {\cal H}_\phi)$ 
on a reproducing kernel Hilbert space ${\cal H}_\phi \subeq {\cal O}(S)$ by 
$(\pi(s).f)(x) = f(xs)$, where the reproducing kernel 
is given by $K(s,t) = \phi(st^*)$ for some holomorphic function 
$\phi \in {\cal H}_\phi$. 
\end{enumerate}
\end{proposition}

\begin{proof} (a) \cite{Ne99}, Prop.~IV.2.3; (b) \cite{Ne99}, Prop.~II.2.11(ii); 
(c) \cite{Ne99}, Lemma~IV.2.6. 
\end{proof}

\section{Host semigroups and host algebras} 

In this section 
we describe the connection between Lie groups and complex semigroups. 
The key point is that there is a Lie theoretic notion of a 
host semigroup $S$ of a Lie group $G$ which can be used to obtain host 
algebras for~$G$. We start with the definition of a 
host semigroup of a Lie group and turn in the second subsection to 
host algebras of topological groups. 

\subsection{Host semigroups of Lie groups} 

\begin{definition}
  \label{def:1.9} 
Let $S$ be a complex involutive
semigroup. A {\it multiplier} of $S$ is a pair $(\lambda, \rho)$ of holomorphic 
mappings $\lambda, \rho \: S \to S$ satisfying the following
conditions:
$$ a \lambda(b) = \rho(a)b, \qquad \lambda(ab) = \lambda(a)b, \qquad 
\mbox{ and } \quad \rho(ab) = a \rho(b).$$

\nin We write $M(S)$ for the set of all multipliers of $S$ and turn
it into an involutive semigroup by 
$$(\lambda, \rho) (\lambda', \rho') := (\lambda \circ \lambda', \rho'
\circ \rho)  \quad \hbox{ and } \quad 
(\lambda, \rho)^* := (\rho^*, \lambda^*), $$
where $\lambda^*(a) := \lambda(a^*)^*$ and 
$\rho^*(a) = \rho(a^*)^*$ (cf.\ \cite{FD88}, p.778). 
\end{definition}

\begin{remark}
  \label{rem:1.10} The assignment 
$\eta_S \: S \to M(S), a \mapsto(\lambda_a, \rho_a)$ 
defines a morphism of involutive semigroups which
is surjective if and only if $S$ has an identity. 
Its image is an involutive semigroup ideal in $M(S)$. 
\end{remark}

\begin{proposition} \label{prop:multi} Let 
$\alpha$ be a locally bounded absolute value 
on the complex involutive semigroup $S$. 
Then the following assertions hold: 
\begin{description}
\item[\rm(1)] For each non-degenerate $\alpha$-bounded holomorphic 
representation $(\pi, {\cal H})$ of $S$ there exists a unique 
unitary representation $(\tilde\pi,{\cal H})$ of $U(M(S))$, determined by 
$\tilde\pi(g)\pi(s) = \pi(gs)$ for $g \in U(M(S))$, $s \in S$. 
\item[\rm(2)] There exists a unique homomorphism 
$\tilde\eta \: U(M(S)) \to U(M(C^*(S,\alpha)))$ with 
$\tilde\eta(g)\eta(s) = \eta(gs)$ for $g \in U(M(S))$, $s \in S$. 
\end{description}
\end{proposition}

\begin{proof} (a) Every $\alpha$-bounded holomorphic representation 
is a direct sum of cyclic ones 
which in turn are of the form $(\pi_\phi, {\cal H}_\phi)$ 
(Proposition~\ref{prop:1.3}). 
We therefore may assume that $\pi = \pi_\phi$ is realized 
on a reproducing kernel space ${\cal H}_\phi \subeq {\cal O}(S)$ 
with reproducing kernel $K(s,t) :=\phi(st^*)$. Then 
$K$ is invariant under the right action of any 
$g = (\lambda_g,\rho_g) \in U(M(S))$: 
$$ K(sg,tg) = \phi((sg)(tg)^*) 
= \phi(sgg^{-1}t^*) = \phi(st^*) = K(s,t). $$
Hence $\tilde\pi_\phi(g)(f) := f \circ \rho_g$ defines a unitary operator 
on ${\cal H}_\phi$ satisfying 
$\tilde\pi_\phi(g)\pi_\phi(s) = \pi_\phi(gs)$ for $s \in S$ 
(cf.\ \cite{Ne99}, Remark~II.4.5). 

(b) Since there exists a faithful representation 
$\pi \: C^*(S,\alpha) \to B({\cal H})$, this follows directly 
from (a). 
\end{proof}

\begin{definition} (a) For a Lie group $G$ 
with Lie algebra $\L(G)$, we call a 
smooth function $\exp_G \: \L(G) \to G$ an 
{\it exponential function} if for each $x \in \L(G)$ the curve 
$\gamma_x(t) := \exp_G(tx)$ is a one-parameter group 
with $\gamma_x'(0) = x$. 
In general an exponential function need not exist, but it is unique 
(\cite{GN07}). 

(b) A Lie group $G$ with an exponential function 
is called {\it locally exponential} if there exists an open  $0$-neighborhood $U$ 
in $\L(G)$ for which $\exp_G\res_U$ is a diffeomorphism onto an open 
subset of $G$. 
\end{definition}

\begin{definition} We say that 
a net $(u_i)_{i \in I}$ in a topological involutive 
semigroup $S$ is an {\it approximate identity}  
if $ \lim u_i s = \lim su_i = s$ holds for all $s \in S$. 
\end{definition}

\begin{remark} For any complex involutive semigroup $S$ 
with an approximate identity, the 
natural map $S \to M(S)$ is injective. In fact, if $(u_i)$ is an 
approximate identity of $S$, then the assertion follows from 
$\eta(s)u_i = su_i \to s$. We may thus identify $S$ with a subsemigroup 
of $M(S)$. 
\end{remark}

\begin{definition} \label{def:hostsem} Let 
$G$ be a connected Lie group with a smooth exponential 
function $\exp_G \: \L(G) \to G$. A triple $(S,\eta,W)$, consisting 
of a complex involutive semigroup $S$ with an approximate identity, a 
group homomorphism $\eta \: G \to U(M(S))$ into the group 
of unitary holomorphic multipliers of $S$ and an open convex 
$\Ad(G)$-invariant cone $W \subeq \L(G)$ 
is called a {\it host semigroup for $G$} 
if the following conditions are satisfied: 
\begin{description}
\item[\rm(HS1)] The left action of $G$ on $S$ defined by $\eta$ is smooth. 
\item[\rm(HS2)] For each $x \in W$, the one-parameter group 
$$\eta_x \: \R \to U(M(S)), \quad t \mapsto \eta(\exp_G(tx)) $$
extends to a morphism 
$$ \hat\eta_x \: \C_+ = \R + i [0,\infty[ \to M(S) $$
of involutive semigroups defining a continuous left action of 
the closed upper halfplane 
$\C_+$ on $S$, $\hat\eta_x(\C_+^0) \subeq S$ (considered as a subsemigroup 
of $M(S)$), and the corresponding map $\C_+^0 \to S$ is holomorphic. 
\item[\rm(HS3)] If $f \: S \to \C$ is a holomorphic function for which all 
functions $f \circ \gamma_x^S$ vanish on the open upper half plane, then 
$f = 0$.  
\end{description}
\end{definition}

\begin{remark} Suppose that $(S,\eta, W)$ is a host semigroup of $G$ 
and $x \in W$. Then the one-parameter subsemigroup 
$\hat\eta_x(it)$, $t > 0$, is an approximate identity for $S$ in the sense that 
for each $s \in S$ we have 
$$ \lim_{t \to 0} \hat\eta_x(it)s = \lim_{t \to 0} s\hat\eta_x(it) = s $$
because the left action of the closed halfplane $\C_+$ on $S$ defined by 
$\hat\eta_x$ is continuous. 
\end{remark}

\begin{proposition} \label{prop:hostexam} 
Let $G_\C$ be a connected complex locally 
exponential Lie group 
whose Lie algebra $\L(G)_\C$ is the complexification of the real Lie algebra 
$\L(G)$, $\sigma$ a holomorphic involutive automorphism 
of $G_\C$ 
with $\L(\sigma)(x + iy) = x-iy$, and $G := (G_\C)^\sigma_0$ the identity 
component of the group $G_\C^\sigma$ of $\sigma$-fixed points in $G_\C$. 
Let $S \subeq G_\C$ be an open connected subsemigroup 
invariant under the involution $s^* := \sigma(s)^{-1}$ 
and $W \subeq \L(G)$ an open convex invariant cone with 
$$ \exp_{G_\C}(iW) \subeq S \quad \mbox{ and } \quad G S G = S.$$ 

Then we obtain for each $g \in G$ a holomorphic multiplier 
$\eta(s) \in U(M(S))$ by 
$\eta(g) = (\lambda_g, \rho_g)$ and $(S,\eta, W)$ is a host semigroup for $G$. 
\end{proposition} 

\begin{proof} (HS1) follows from the smoothness of the multiplication in $G_\C$. 

\nin (HS2): For $x \in W$ we put 
$\hat\eta_x(z) := (\lambda_{\exp(zx)}, \rho_{\exp zx})$ and 
note that this is the multiplier corresponding to the semigroup element 
$\exp(zx)$ because for $z = a + i b$ we have 
$\exp(zx) = \exp(ax)\exp(ibx)\in G S \subeq S$. 

\nin (HS3): Let $f \: S \to \C$ be a holomorphic function vanishing 
on the sets $\exp(\C_+^0 x)$, $x \in W$. 
Let $\Omega \subeq \exp^{-1}(S) \subeq \L(G)_\C$ 
denote the connected component containing 
the cone $i W$. Then the holomorphic function $f \circ \exp \: 
\Omega\to \C$ (cf.\ \cite{GN07} for the holomorphy of $\exp$) vanishes 
on $iW$, hence on a neighborhood of $iW$, and therefore on all of $\Omega$. 
Since $G_\C$ is locally exponential, there exists a point 
$x_0 \in i W$ (sufficiently close to~$0$) and an open neighborhood 
$U$ of $x_0$ in $\Omega$, such that $\exp(U)$ is an open subset of $S$. 
Then $f$ vanishes on $\exp(U)$ and 
hence on all of $S$ because $S$ is connected. 
\end{proof}

\begin{example} \label{ex:tube} 
Let $G = V$ be a locally convex space, $W \subeq V$ an open 
convex cone and $S := V + i W$. Then $\eta(v)(s) := s + v$ yields 
a host semigroup $(S,\eta,W)$ for $V$. 
\end{example}

\subsection{Host algebras} 

\begin{definition} If $\cA$ is a $C^*$-algebra, then we write 
$M(\cA)$ for the set of continuous linear multipliers on $\cA$. 
Then $M(\cA)$ carries a natural structure of a $C^*$-algebra and 
the map $\eta_\cA \: \cA \to M(\cA)$ is injective 
(cf.\ \cite{Pe79}, Sect.~3.12). 
We write ${\cal A}$ for its
image in $M({\cal A})$. The  {\it strict topology} on $M(\cA)$  
is the locally convex topology defined by the seminorms 
$$ p_a(m) := \|m a\| + \|a m\|,
\qquad a\in {\cal A}, m\in M({\cal A}).$$ 
The involution is continuous with respect to this topology and 
the multiplication is continuous on bounded subsets, which implies 
in particular that the unitary group $U(M(\cA))$ is a topological group 
(cf.\ \cite{Wo95}, Sect.~2). 

For a complex Hilbert space ${\cal H}$, we write $\Rep({\cal A},{\cal H})$ 
for the set of non-degenerate representations of ${\cal A}$ on ${\cal H}$. 
Each representation $(\pi, {\cal H})$ of ${\cal A}$ which is 
non-degenerate in the sense that $\pi(\cA)v = \{0\}$ implies $v = 0$ 
extends to a unique representation $\tilde\pi$ of 
$M(\cA)$ satisfying $\tilde\pi \circ \eta_\cA = \pi$ 
which is continuous with respect to the 
{\it strict topology} on $M({\cal A})$ and 
the strong operator topology on $B({\cal H})$ 
(cf.\ Proposition~\ref{prop:b.1} below). 
\end{definition}

\begin{examples} \label{ex:multi} (\cite{Pe79}, Sect.~3.12)  
(a) If $\cA$ is a closed $*$-subalgebra of $B({\cal H})$, 
then $M(\cA) \cong \{ X \in B({\cal H}) \: X \cA + \cA X \subeq \cA\}$ 

\nin(b) If $\cA = K({\cal H})$ is the ideal of compact operators in 
$B({\cal H})$, then $M(K({\cal H})) \cong B({\cal H})$. 

\nin(c) If $\cA = C_0(X)$ is the $C^*$-algebra of continuous functions vanishing 
at infinity on the locally compact space $X$, then 
$M(\cA) \cong C_b(X)$ is the $C^*$-algebra of bounded continuous functions 
on $X$. 
\end{examples}

\begin{definition} Let $G$ be a topological group. 
A {\it host algebra for $G$} is a pair 
$({\cal A}, \eta)$, where  ${\cal A}$ is a $C^*$-algebra and 
$\eta \: G \to U(M({\cal A}))$ is a group homomorphism 
such that: 
\begin{description}
\item[\rm(H1)] For each non-degenerate representation $(\pi, {\cal H})$ 
of $\cA$, the representation $\tilde\pi \circ \eta$ of $G$ is continuous. 
\item[\rm(H2)] For each complex Hilbert space 
${\cal H}$, the corresponding map 
$$ \eta^* \: \Rep({\cal A},{\cal H}) \to \Rep(G, {\cal H}), \quad 
\pi \mapsto \tilde\pi \circ \eta $$ 
is injective. 
\end{description}
\nin We say that $(A,\eta)$ is a 
{\it full host algebra} if 
$\eta^*$ is surjective for each Hilbert space~${\cal H}$.  
\end{definition}

\begin{remark} \label{rem:2.2} 
(a) If $\eta \: G\to U(M(\cA))$ is {\it strictly continuous}, i.e., 
continuous with respect to the strict topology on $U(M(\cA))$, then 
Proposition~\ref{prop:b.1}(3) implies (H1). 

(b) Since the extension $\tilde\pi$ of a non-degenerate representation 
$\pi$ of $\cA$ is strictly continuous, condition (H2) holds 
if $\eta(G)$ spans a strictly dense subalgebra of $M(\cA)$. 
In view of \cite{Wo95}, Prop.~2.2, (H2) conversely implies that 
$\Spann(\eta(G))$ is strictly dense in $M(\cA)$. 

(c) For any multiplier action $\eta \: G\to U(M(\cA))$ of a 
topological group $G$ on the $C^*$-algebra $\cA$, the subspace 
$R$, consisting of all elements $a \in \cA$ for which the map 
$G \mapsto \cA, g \mapsto \eta(g)a$ is continuous is a right ideal 
which is closed because $G$ acts by isometries on $\cA$. It is biinvariant 
under $G$. 
Hence $\cA_c := R \cap R^*$ is a $C^*$-subalgebra of $\cA$ which is 
$G$-biinvariant, and the corresponding homomorphism 
$\eta_c \: G \to U(M(\cA_c))$ is strictly continuous. 

(d) If a homomorphism 
$\eta \: G \to U(M(\cA))$ satisfies (H2), then (b) implies 
that $\eta(G)$ spans a strictly dense subalgebra of $M(\cA)$. 
This implies in particular, that the $G$-biinvariant closed subalgebra 
$\cA_c$ of $\cA$ is a two-sided ideal of $M(\cA)$ and the corresponding 
morphism $\gamma \: M(\cA) \to M(\cA_c)$ is obviously 
strictly continuous and satisfies $\gamma \circ \eta = \eta_c$. 

We claim that $(\cA_c, \eta_c)$ 
also is a host algebra for $G$. In fact, (H1) follows from the strict 
continuity of $\eta_c$. Next we note that 
$\gamma(M(\cA))$ contains $\cA_c$, so that it is strictly 
dense in $M(\cA_c)$ (\cite{Wo95}, Prop.~2.2). Therefore 
$\eta_c(G) =  \gamma(\eta(G))$ spans a strictly dense subalgebra of 
$M(\cA_c)$, which implies (H2), as we have seen in (b). 
\end{remark}

\begin{examples} \label{exs:hosts} (a) Let 
$G$ be a locally compact group and $C^*(G)$ the enveloping 
$C^*$-algebra of the group algebra $L^1(G)$. Then we 
have a natural homomorphism 
$\eta \: G \to U(M(C^*(G)))$ which is determined by the left action 
$\eta(g)(f)(x) = f(g^{-1}x)$ on $L^1$-functions. 
Since $G$ acts continuously from the left and 
the right on $L^1(G)$ and the image of $L^1(G)$ is dense in 
$C^*(G)$, $\eta$ is continuous with respect to the strict topology. 
It is well known that $(C^*(G),\eta)$ is a full host algebra of 
$G$ (\cite{Dix64}, Sect.~13.9). 

(b) Let $G$ be an abelian topological group and $\hat G := \Hom(G,\T)$ its 
character group. Then any host algebra 
$(\cA,\eta)$ for $G$ is commutative because of the strict density of 
$\eta(G)$ in $M(\cA)$. Hence there exists a locally compact space 
$X$ with $\cA \cong C_0(X)$ and $M(\cA) \cong C_b(X)$ 
(cf.\ Example~\ref{ex:multi}). Then $U(M(\cA)) \cong C(X,\T)$, 
where the strict topology on this group corresponds to the compact  
open topology and a $*$-subalgebra of $C_b(X)$ is strictly 
dense if and only if it separates the points of $X$ 
(\cite{Br77}, Lemma~3.5). Therefore the map 
$\gamma \: X \to \hat G$ defined by 
$\gamma(x)(g) := \eta(g)(x)$ is injective, 
so that we may consider $X$ as a subset of the character group $\hat G$. 

If, conversely, $X \subeq \hat G$ is a subset, endowed with a locally 
compact topology finer than the topology of pointwise convergence on $G$, 
then the natural map $\eta \: G \to C(X,\T) = U(M(C_0(X)))$ defined by 
$\eta(g)(\chi) := \chi(g)$ satisfies (H2) because $\eta(G)$ 
separates the points of $X$. If, in addition, $\eta$ is strictly continuous, 
i.e., each compact subset of $X$ is equicontinuous, then (H1) is 
also satisfied, so that $(C_0(X),\eta)$ is a host algebra of $G$. 

(c) If $\cA$ is any $C^*$-algebra and $G := U(M(\cA))$, endowed with 
the strict topology, then the fact that $G$ spans $M(\cA)$ implies that 
$\eta = \id_G$ satisfies (H1) and (H2), so that $(\cA,\id_G)$ 
is a host algebra for $G$. 
\end{examples} 

\begin{example} (Non-uniqueness of host algebras) 
Let $G := \Z$. Then its character group is $\hat G \cong \T$, which is a 
compact group with respect to the topology of pointwise convergence.
Since $G$ is locally compact, $C^*(G) \cong C(\T)$ is a full host algebra 
for $G$. Let $\cA := C_0([0,1[)$ and define a homomorphism 
$\eta \: \Z \to U(M(C_0([0,1[))) \cong C([0,1[,\T)$ by 
$\eta(n)(x) := e^{2\pi i nx}$. Then $\eta(1) \: [0,1[ \to \T$ is a 
continuous bijection, which implies in particular that 
$\eta(\Z)$ separates the points, so that (H2) holds. 
Further, $\Z$ is discrete, so that (H1) is trivially satisfied, 
and thus $(\cA,\eta)$ is a host algebra. 
This host algebra is full because the representations of $\Z$ are in one-to-one correspondence with 
Borel spectral measures on $\T$ and $\eta(1)$ is a Borel isomorphism. 
Note in particular that the full host algebra $\cA$ is not unital, although 
$G$ is a discrete group. 
\end{example}

\begin{remark}  Let $G$ be a Lie group, $\cA$ a $C^*$-algebra and 
$\eta \: G \to U(M(\cA))$ a group homomorphism. We then obtain 
left and right actions of $G$ on $\cA$ by isometries. 
Let $\cA^\infty \subeq \cA$ denote the set of all elements for which the 
orbit maps of these actions are smooth. 
The space $\cA^\infty_l$ of smooth vectors for 
the left action $\eta_l$ is a right ideal, the set $\cA^\infty_r$ of smooth 
vectors for the right action $\eta_r$ is a left ideal and both 
are exchanged by the involution. Hence their intersection $\cA^\infty$ is a 
$*$-subalgebra on which $G$ acts by multipliers. 
\end{remark}

\begin{definition} \label{def:2.x} 
Let $G$ be a Lie group and $\cA$ a $C^*$-algebra. 
We say that a homomorphism 
$\gamma \: G \to U(M(\cA))$ is {\it strictly smooth}  
if the $*$-subalgebra $\cA^\infty$ of smooth vectors for the left and right 
action of $G$ on $\cA$ is dense. 

If this is the case, then the maximal $C^*$-subalgebra $\cA_c$ of $\cA$ 
on which $G$ acts continuously is dense, hence coincides with $\cA$, 
so that $\eta$ is in particular strictly continuous 
(cf.\ Remark~\ref{rem:2.2}(c)). 
\end{definition}

\begin{proposition} \label{prop:holhostalg} 
Let $(S,\eta,W)$ be a host semigroup 
of the Lie group $G$ and $\alpha$ a $G$-invariant locally bounded 
absolute value on $S$. Then $\eta$ induces a strictly smooth homomorphism 
$\tilde\eta \: G \to U(M(C^*(S,\alpha)))$ determined uniquely by 
$\tilde\eta(g)(\eta_\alpha(s)) = \eta_\alpha(gs)$ for $g \in G$, $s\in S$, 
and $(C^*(S,\alpha),\tilde\eta)$ is a host algebra for $G$. 
\end{proposition}

\begin{proof} The existence of $\tilde\eta$ follows 
from Proposition~\ref{prop:multi}, since $U(M(S))$ acts by 
unitary multipliers on $C^*(S,\alpha)$. That $\tilde\eta$ defines 
a strictly smooth multiplier action of $G$ on $C^*(S,\alpha)$ 
follows from the relation $\tilde\eta(g)\eta_\alpha(s) = \eta_\alpha(gs)$, 
which implies that $\eta_\alpha(S)$ consists of smooth vectors for $G$ because 
$\eta_\alpha \: S \to C^*(S,\alpha)$ is a holomorphic map. 
Hence $\tilde\eta$ is strictly smooth. 

To see that $\tilde\eta$ defines a host algebra, let 
$(\pi_i, {\cal H})$, $i = 1,2$, be two representations of 
$C^*(S,\alpha)$ with 
$\tilde\pi_1 \circ \tilde\eta = \tilde\pi_2 \circ \tilde\eta$. 
For each $x \in W$ we then have 
$$ \tilde\pi_1 \circ \tilde\eta_x = \tilde\pi_2 \circ \tilde\eta_x 
\: \R \to U({\cal H}).$$
Now $\tilde\pi_i \circ \eta_x \: \C_+ \to B({\cal H})$, $i =1,2$,  
are two continuous representations of the involutive semigroup $\C_+$ which 
are holomorphic on $\C_+^0$ and coincide on $\R$, so that 
\cite{Ne99}, Lemma~XI.2.2, implies that they are equal. 
Hence $\pi_1 \circ \eta_x = \pi_2 \circ \eta_x$ for each $x \in W$, so 
that $\pi_1 - \pi_2 \: S \to B({\cal H})$ is a holomorphic function 
vanishing on all sets $\gamma_x^S(\C_+^0)$, and now (HS3) leads to 
$\pi_1 = \pi_2$. 
\end{proof}

\section{Multiplier actions of Lie groups on \break $C^*$-algebras} 

In the preceding section we have seen that we can associate to each 
host semigroup $S$ and any $G$-invariant locally bounded absolute 
value $\alpha$ on $S$ a host algebra $C^*(S,\alpha)$ for $G$. 
In this section we slightly change our perspective 
and ask for properties of a homomorphism $\eta \: G \to U(M(\cA))$ 
which are characteristic for a host algebras of the form $\cA = C^*(S,\alpha)$. 
This will lead us to the momentum set $I_\eta \subeq \L(G)'$ 
of the pair $(\cA,\eta)$. We shall see in particular that the 
weak-$*$-closed convex $\Ad^*(G)$-invariant set $I_\eta$ tells us 
for which open invariant cones $W \subeq \L(G)$ there might be 
a corresponding host semigroup. 

\subsection{Strictly continuous multiplier actions} 

If $G$ is a topological group and $\cA$ a $C^*$-algebra, we also call 
a strictly continuous homomorphism $\eta \: G \to U(M(\cA))$ a 
{\it strictly continuous multiplier action of $G$ on $\cA$}. 

\begin{remark}   \label{rem:3.1} 
(a) If $G$ is a locally compact group, $\cA$ is a $C^*$-algebra and the 
homomorphism $\gamma \: G \to U(M(\cA))$ is strictly continuous, then 
integration yields a morphism 
$$ \gamma \: L^1(G) \to M(\cA), \quad \gamma(f) := 
\int_G f(g)\gamma(g)\ dg $$
of Banach-$*$-algebras, 
so that the universal property of the $C^*$-algebra $C^*(G)$, 
the enveloping $C^*$-algebra of $L^1(G)$, implies the 
existence of a corresponding morphism of $C^*$-algebras 
$\tilde\gamma \: C^*(G) \to M(\cA)$
for which $\tilde\gamma(C^*(G))\cA \supeq \gamma(L^1(G))\cA$ 
is dense in $\cA$ 
(use an approximate identity in $L^1(G)$ and the 
strict continuity of the action of $G$)
and we have $\tilde\gamma \circ \eta = \gamma.$

If, conversely, $\alpha \: C^*(G) \to M(\cA)$ is a morphism of 
$C^*$-algebras for which $\alpha(C^*(G))\cA$ is dense in $\cA$, 
then Proposition~\ref{prop:b.2} in the appendix implies that 
$\alpha$ extends to strictly continuous morphism 
$\tilde\alpha \: M(C^*(G)) \to M(\cA)$. Hence 
$\gamma := \tilde\alpha \circ \eta \: G \to U(M(\cA))$ is a strictly continuous 
homomorphism with $\tilde\gamma = \tilde\alpha$. 

We conclude that strictly continuous multiplier actions of $G$ on a 
$C^*$-algebra $\cA$ are in one-to-one correspondence with 
morphisms $\alpha \: C^*(G) \to M(\cA)$ for which 
$\alpha(C^*(G))\cA$ is dense in $\cA$. 
\end{remark}

\begin{example} Let $\cA := K({\cal H})$ be the $C^*$-algebra of 
compact operators on the complex Hilbert space ${\cal H}$. 
Then $M(\cA) \cong B({\cal H})$ is the algebra of all bounded operators on 
${\cal H}$. Hence $U(M(\cA)) \cong U({\cal H})$ is the unitary group 
of ${\cal H}$ and the strict topology on this group coincides 
with the strong operator topology because 
if $U({\cal H})$ carries the strong operator topology, 
the closed subalgebra $K({\cal H})_c$ contains 
finite rank operators, hence coincides with $K({\cal H})$. 
Therefore a strictly continuous multiplier action of a topological 
group $G$ on $K({\cal H})$ is the same as a continuous unitary representation 
on ${\cal H}$. 
\end{example}

\begin{remark} If $G$ is a finite dimensional Lie group 
and $\eta \: G \to U(M(\cA))$ 
is strictly continuous, then it is also strictly smooth. In fact, 
there exists a sequence $(\delta_n)_{n \in \N} \in C_c^\infty(G,\R)$ 
which is an approximate identity in $L^1(G)$. 
Then each $a \in \cA$ is the norm limit of the elements 
$\eta(\delta_n)a\eta(\delta_n)\in \cA^\infty$. 
\end{remark}

\subsection*{Momentum sets of unitary representations} 

\begin{definition} Let $(\pi, {\cal H})$ be a continuous unitary
representation of the Lie group $G$ on the Hilbert space ${\cal H}$. 
We write ${\cal H}^\infty \subeq {\cal H}$ for the subspace of smooth 
vectors and note that this is a linear subspace on which we have a 
derived representation $d\pi$ of the Lie algebra $\g = \L(G)$. 
We call the representation 
$(\pi, {\cal H})$ {\it smooth} if ${\cal H}^\infty$ is dense in ${\cal H}$.

\par (a) Let $\bP({\cal H}^\infty) = \{ [v] := 
\C v \: 0 \not= v \in 
{\cal H}^\infty\}$ 
denote the projective space of the subspace ${\cal H}^\infty$ 
of smooth vectors. The map 
$$ \Phi_\pi \: \bP({\cal H}^\infty)\to \g' \quad \hbox{ with } \quad 
\Phi_\pi([v])(x) 
= \frac{1}{i}  \frac{\la  d\pi(x).v, v \ra}{\la v, v \ra} $$
is called the {\it momentum map of the unitary representation $\pi$}.  
The right hand side is well defined because it
only depends on $[v] = \C v$. The operator $i\cdot d\pi(x)$ is symmetric so
that the right hand side is real, and since $v$ is a smooth vector, 
it defines a continuous linear functional on $\g$. 

\par (b) The weak-$*$-closed convex hull 
$I_\pi \subeq \g'$ of the image of $\Phi_\pi$ is called the 
{\it convex momentum set of $\pi$}. 
\end{definition}

\begin{remark} \label{rem:3.5} If $\pi \: \cA \to B({\cal H})$ is a non-degenerate 
representation of a $C^*$-algebra and $\eta \: G \to U(M(\cA))$ a 
strictly smooth multiplier action, 
then the corresponding representation $\tilde\pi \circ \eta$ 
of $G$ on ${\cal H}$ has a dense space ${\cal H}^\infty$ of smooth 
vectors because the dense subspace spanned by 
$\pi(\cA^\infty){\cal H}$ consists of smooth vectors for $G$. 
\end{remark}

\newcommand{\derat}[1]{ \hbox{$\frac{d}{dt}$\vrule}_{t=#1}}  

\begin{lemma} \label{lem:reed-simon} Let 
$G$ be a Lie group with exponential function 
and $(\pi, {\cal H})$ a unitary representation with a dense space 
${\cal H}^\infty$ of smooth vectors. Then, for each $x \in \L(G)$, 
the unbounded operator 
$$ d\pi(x)  \: {\cal H}^\infty \to {\cal H}, \quad 
d\pi(x)v := \derat0 \pi(\exp_G(tx)).v $$
is essentially skewadjoint and its closure $\oline{d\pi(x)}$ 
is the infinitesimal generator of the unitary one-parameter group 
$\pi_x := \pi \circ \gamma_x$, $\gamma_x(t) = \exp_G(tx)$. 
\end{lemma} 

\begin{proof} Since the dense subspace ${\cal H}^\infty$ is invariant 
under $\pi(\gamma_x(\R))$, the assertion follows from 
\cite{RS75}, Thm.~VIII.10. 
\end{proof}

\begin{lemma} \label{lem:spec-onepar} If $(\pi, {\cal H})$ 
is a smooth unitary representation of the Lie group $G$
with exponential function, $x \in \L(G)$, 
$\pi_x := \pi \circ \gamma_x$ and $A_x := -i \oline{d\pi_x(1)}$ 
the corresponding selfadjoint operator, then 
$$ \inf \Spec(A_x) = \inf \la I_\pi, x \ra. $$
\end{lemma}

\begin{proof} For $m(x) 
:= \inf \la I_\pi, x\ra \in \R \cup \{-\infty\}$, we have 
$$ \la A_x.v,v\ra \geq m(x) \la v, v \ra \quad \mbox{ for each } \quad 
v \in {\cal H}^\infty, $$
and since the graph 
of the operator $-i \cdot d\pi(x)$ 
on ${\cal H}^\infty$ in dense in the graph of 
$A_x$ (Lemma~\ref{lem:reed-simon}), 
$\la A_x.v,v\ra \geq m(x) \la v, v \ra$ holds 
for each $v$ in the domain of $A_x$. This shows that 
$\inf \Spec(A_x) \geq m(x)$. The converse inequality holds trivially. 
\end{proof}

\begin{problem} Let $(\pi, {\cal H})$ be a unitary representation of 
the Lie group $G$ on ${\cal H}$. If $v$ is a smooth vector, then the 
function $\pi^v(g) := \la g.v, v \ra$ is smooth. In 
\cite{Ne99}, Prop.~X.6.4 it is shown that the converse also holds if 
$G$ is finite dimensional. Does this result generalize to infinite dimensional 
Lie groups? 
\end{problem}

\subsection*{Holomorphic extension of multiplier actions} 

\begin{definition} Let $\eta \: G \to U(M(\cA))$ be a strictly smooth 
multiplier action of $G$ on the $C^*$-algebra $\cA$ and $\g = \L(G)$. 

We write $S(\cA)$ for the 
set of states of $\cA$. Since each state of $\cA$ is of the form 
$\phi(a) = \pi^v(a) := \la \pi(a).v,v\ra$ for a unit vector $v \in {\cal H}$ 
and a non-degenerate representation $(\pi, {\cal H})$ of $\cA$, 
there exists a canonical extension 
$\tilde\phi := \tilde\pi^v$ to a state of $M(\cA)$. 

We call a state $\phi$ of $\cA$ {\it $\eta$-smooth} if 
$\tilde\phi \circ\eta$ is smooth 
and write $S(\cA)^\infty$ for the set of 
$\eta$-smooth states of $\cA$. We now have a 
{\it momentum map} 
$$ \Psi_\eta \: S(\cA)^\infty \to \g', \quad 
\Psi_\eta(\phi) = \frac{1}{i} d(\tilde\phi \circ \eta)(\1). $$
Since $S(\cA)^\infty$ is a convex set, the weak-$*$-closure 
$$I_\eta := \oline{\Psi_\eta(S(\cA)^\infty)} \subeq \g' $$
also is a convex subset of $\g'$; called the momentum set of $(\cA,\eta)$.
\end{definition}

\begin{proposition} \label{prop:holext} 
Let $\eta \: G \to U(M(\cA))$ be a strictly smooth multiplier 
representation, $I_\eta \subeq \g'$ its momentum set and $m \in \R$. 
Then the following are equivalent for $x \in \g$: 
\begin{enumerate}
\item[\rm(1)] If $\eta_x(t) := \eta(\exp_G(tx))$, then 
the corresponding homomorphism of $C^*$-algebras 
$\tilde\eta_x \: C^*(\R) \cong C_0(\R) 
\to M(\cA)$ factors through the quotient algebra $C_0([m,\infty[)$. 
\item[\rm(2)] $I_\eta(x) \geq m$. 
\item[\rm(3)] $\eta_x$ extends to a strictly continuous homomorphism 
$\hat\eta_x \: \C_+ \to M(\cA)$ of involutive semigroups 
which is holomorphic on $\C_+^0$ and satisfies \break 
$\|\hat\eta_x(z)\| \leq e^{-m \Im z}$. 
\end{enumerate}
If these conditions are satisfied, then 
$\|\hat\eta_x(a+ib)\| = e^{-b \cdot \inf I_\eta(x)}$. 
\end{proposition}

\begin{proof} (1) $\Leftrightarrow$ (2): 
Let $(\pi, {\cal H})$ be a universal representation 
of $\cA$, i.e., each 
state of $\cA$ is of the form $\pi^v(a) = \la \pi(a).v,v\ra$ 
for some unit vector 
$v \in {\cal H}$. Then $\tilde\pi \circ \eta$ is a smooth representation of 
$G$ (Remark~\ref{rem:3.5}) 
and $\eta_x$ also defines a continuous unitary representation 
$\pi_x := \tilde\pi \circ \eta_x$ of $\R$ on ${\cal H}$. 

For any smooth unit vector $v \in {\cal H}^\infty$ and the 
corresponding smooth state $\pi^v$ we then obtain with 
$(\tilde\pi^v\circ \eta)(\exp_G(tx)) = \pi_x(t)$: 
$$ \Psi_\eta(\pi^v)(x) = - i \cdot d(\tilde \pi^v \circ \eta)(x) 
= - i \la d\pi(x).v, v \ra = \Phi_\pi([v])(x), $$
which leads to $I_\pi(x) \subeq I_\eta(x)$. 
In view of \cite{Ne99}, Prop.~X.6.4, the smoothness of a state 
$\pi^v$ implies the smoothness of $v$ for $\pi_x$, so that 
$I_\eta(x) \subeq I_{\pi_x}$. 
Lemma~\ref{lem:spec-onepar} now implies that 
$\inf I_{\pi_x} = \inf \la I_\pi, x \ra = \inf I_\pi(x)$, so that we 
arrive at 
$$ \inf I_\pi(x) = \inf I_\eta(x) = \inf I_{\pi_x}.$$

A simple argument with the spectral measure of the unitary 
one-parameter group $\pi_x$ shows that the 
kernel of the corresponding representation $\hat\pi_x$ of 
$C^*(\R) \cong C_0(\R)$ contains the ideal 
$$I_m := \{ f \in C_0(\R) \: \supp(f) \subeq ]-\infty,m[\}$$ 
if and only if $m \leq I_{\pi_x}$, which means that it 
factors through the quotient algebra $C_0(\R)/I_m \cong C_0([m,\infty[).$
If this is the case, then the image of $\hat\pi_x$ 
lies in the multiplier algebra 
$M(\cA) \cong \{ T \in B({\cal H}) \: T \cA + \cA T \subeq \cA \}$
(Proposition~\ref{prop:b.2}). This proves the equivalence of (1) and (2). 

(1) $\Rarrow$ (3): First we consider the map 
$$ \gamma \: \C^+ \to C_b([m,\infty[) = M(C_0([m,\infty[)), \quad 
\gamma(z)(t) := e^{izt}. $$
Then $\|\gamma(z)\| = e^{-m\Im z}$, so that $\gamma$ is locally bounded. 
Since the strict topology 
on bounded subsets of $C_b([m,\infty[)$ coincides with the 
compact open topology (cf.\ \cite{Br77}, Lemma~3.5), 
the strict continuity 
of $\gamma$ follows from the continuity of the 
map $\beta \: \C \to C_b([m,\infty[),\beta(z)(t) := e^{izt}$ 
with respect to the compact open topology. That 
$\gamma$ is holomorphic on the open upper halfplane $\C_+^0$ 
is a consequence of Example~\ref{ex:cstarsem}(c). Now (3) 
follows by composing the strictly continuous extension 
$M(C_0([m,\infty[)) \cong C_b([m,\infty[) \to M(\cA)$ with~$\gamma$. 

(3) $\Rarrow$ (2): By definition, $\hat\eta_x$ induces a morphism 
$\beta \: C^*(\C_+^0,\alpha) \to M(\cA)$, where 
$\alpha(z) := e^{-m \Im z}$. Since $\hat\eta_x$ is strictly continuous,   
$\beta(\C_+^0)\cA$ is dense in $\cA$, so that $\beta$ extends to a 
strictly continuous morphism 
$\hat\beta \: M(C^*(\C_+^0,\alpha)) \cong C_b([m,\infty[) \to M(\cA)$ 
with $\hat\beta(\gamma(t)) = \eta_x(t)$ for $t \in \R$. Therefore 
the homomorphism $\tilde\eta_x \: C^*(\R) \to M(\cA)$ factors through 
the quotient $C_b([m,\infty[)$. 
\end{proof}

\begin{remark} \label{rem:2.23} 
(a) As the example $\cA = C_0([m,\infty[)$ shows, the map 
$\hat\eta_x$ need not be norm continuous because the natural map 
$$ \gamma \: \C_+ \to C_b([m,\infty[), \quad \gamma(z)(t) = e^{itz}$$ 
is not norm continuous at the boundary $\R = \partial \C_+$. 

(b) Assume that the conditions of Proposition~\ref{prop:holext} are 
for the element $x \in \g$. 
Let $\cB := M(\cA)_c$ denote the $C^*$-subalgebra consisting of all 
elements on which $G$ acts continuously by multipliers from the left and 
the right. Then $\eta_x(\R)\cB + \cB\eta_x(\R) \subeq \cB$ implies that 
$\eta_x$ induces a strictly continuous morphism 
\break $\eta_x^\cB \: \R\to U(M(\cB))$.

Since the induced homomorphism $C^*(\R) \to M(\cA)$ factors through 
\break $C_0([m,\infty[)$, the same holds for the corresponding morphism 
$C^*(\R) \to M(\cB)$. From that we conclude that we even obtain a 
strictly continuous morphism 
$\C_+ \to M(\cB)$ which is holomorphic on $\C_+^0$. 
\end{remark}

\begin{proposition} \label{prop:host-cone} 
Let $\eta \: G \to U(M(\cA))$ be the strictly smooth 
multiplier action defined by a host algebra of $G$ obtained from a 
host semigroup $(S,\eta, W)$ for which the map 
$$ \Exp \: W \to S, \quad x \mapsto \hat\eta_x(i) $$
is continuous. Then 
$$ s \: W \to \R, \quad s(x) := - \inf \la I_\eta, x \ra $$
is a locally bounded function on $W$. 
\end{proposition} 

\begin{proof} For each $x \in W$, the 
homomorphism $\eta_x \: \R \to U(M(C^*(S,\alpha)))$ 
extends to a homomorphism of involutive semigroups 
$\hat\eta_x \: \C_+ \to M(C^*(S,\alpha))$ which is 
holomorphic on $\C_+^0$ and comes from a smooth multiplier action of 
$\C_+$ on $S$ (Definition~\ref{def:hostsem}). We also have 
$$ \|\hat\eta_x(z)\| \leq \alpha(\gamma_x(z)) \quad \mbox{ for } 
\quad z \in \C_+^0. $$
From Example~\ref{ex:cstarsem}(c) we now derive that 
$\|\hat\eta_x(z)\| = e^{-m_x\Im z}$ for some $m_x \in \R$ 
and in particular that $\hat\eta_x$ is locally bounded on $\C_+$. 
Therefore the continuity of the corresponding 
multiplier action of $\C_+$ on $S$ implies the strict continuity of the 
corresponding multiplier action on $C^*(S,\alpha)$. Now 
Proposition~\ref{prop:holext}(3) tells us that $\inf I_\eta(x) =  m_x$, 
so that $s(x) = -m_x$. Further, 
$$ e^{s(x)} = e^{-m_x} 
= \|\hat\eta_x(i)\| = \|\eta_\alpha(\hat\eta_x(i))\| 
\leq \alpha(\Exp(x)) 
$$
shows that $s$ is locally bounded on the open cone $W$. 
\end{proof}

\section{Convex functions on infinite dimensional domains} 

In Proposition~\ref{prop:host-cone} we have seen how host 
algebras of a Lie group $G$ coming from 
host semigroups lead to weak-$*$-closed convex subsets 
$I_\eta$ of the dual $\L(G)'$ of the locally convex Lie algebra $\L(G)$ 
with the property that the support function $x \mapsto - \inf \la I_\eta,x \ra$ 
is locally bounded on some open convex cone in $\L(G)$.

In this section we therefore take a closer look at weak-$*$-closed convex 
subsets $C$ of the dual $V'$ of a locally convex space $V$. 
We are in particular 
interested in conditions for the cone $B(C) = \{ v \in V \: 
\inf \la C, v\ra > -\infty\}$ to have non-empty interior and 
the corresponding support function 
$s_C(v) := - \inf \la C, v\ra$ to be locally bounded on $B(C)^0$. 
We start with a general discussion of convex sets and then turn to 
convex functions in the second subsection. 

\subsection*{Convex subsets of locally convex spaces}

\begin{proposition} \label{prop:3.1}
Let $\eset \not= C \subeq V$ be a closed convex set in the topological vector 
space~$V$.  
\begin{description}
\item[(1)] $\lim(C) :=\{v\in V\colon C+ v\subeq C\}$ is a closed convex cone, 
called the recession cone of $C$. 
\item[(2)] $\lim(C) = \{ v \in V \: v = \lim_{n \to \infty} t_n c_n,
c_n \in C, t_n \to 0, t_n \geq 0 \}$. 
\item[(3)] If $c \in C$ and $x \in V$ satisfy 
$c + \R^+ x \subeq C$, then $x \in \lim(C)$. 
\item[(4)] If $C$ is bounded, then $\lim(C) = \{0\}$. 
\end{description}
\end{proposition}

\begin{proof} (1) The closedness of $\lim(C)$ is an immediate consequence of
the closedness of $C$. 

 (2) If $c \in C$ and $x \in \lim(C)$, then $c + n x \in C$
for $n \in \N$ and  $\frac{1}{n}(c + nx)\to~x. $ 

If, conversely, $x = \lim_{n \to \infty} t_n c_n$ with $t_n \to 0$,
$t_n \geq 0$ and $c, c_n \in C$, then 
$(1 - t_n) c + t_n c_n \to c + x \in \oline C = C$
implies that $C + x \subeq C$, i.e.\ $x \in \lim(C)$. 

(3) In view of (2), this follows from 
$ {1\over n}(c + n x) \to x$. 

(4) If $C$ is bounded, each continuous linear 
functional $f \: V \to \R$ is bounded on $C$. For each $x \in \lim(C)$ 
the relation $C + \N x \subeq C$ then leads to $f(x) = 0$. 
Since $V'$ separates the points of $V$, we obtain $x = 0$. 
\end{proof}

\begin{definition} Let $V$ be a locally convex space and 
$C \subeq V'$ a subset. We put 
$$ B(C) := \{ v \in V \: \inf \la C, v \ra  > - \infty \} 
\quad \hbox{ and } \quad 
C^\star := \{ v \in V \: \la C, v \ra \subeq \R_+\}. $$
Then $C^\star \subeq B(C)$ are convex cones and $C^\star$ is called the 
{\it dual cone of $C$}. If $C$ is a cone, then 
$B(C) = C^\star$. 
\end{definition}

Let $C \subeq V'$ be a weak-$*$-closed convex subset. 
As a consequence of the Hahn--Banach Separation Theorem, there exists 
for each element $\alpha \in V' \setminus C$ some $x \in V$ (the dual of 
$V'$ endowed with the weak-$*$-topology) with \break 
$\alpha(v)  < \inf \la C, v\ra$. Then $v \in B(C)$, and we thus obtain 
\begin{equation}
  \label{eq:reco}
C = \{ \alpha \in V' \: (\forall v \in B(C))\ \alpha(v) \geq \inf \la C, v \ra\},  
\end{equation}
which permits us to reconstruct $C$ from its support function 
$s_C(v) = \break -\inf \la C, v \ra$ on~$V$. 

\begin{lemma} \label{lem:bcdual} 
For a non-empty weak-$*$-closed convex subset $C \subeq V'$, the following 
assertions hold: 
\begin{description}
\item[\rm(i)]  $B(C)$ is a convex cone satisfying $B(C)^\star = \lim(C).$ 
\item[\rm(ii)] If $B(C)$ has non-empty interior,  
then $B(C)$ has the same interior as $\lim(C)^\star$. 
\item[\rm(iii)] If $C$ is a cone, then $B(C) = C^\star$ has non-empty interior 
if and only if $C$ has a weak-$*$-compact equicontinuous base. 
\end{description}
\end{lemma}

\begin{proof} (i) The relation $C + \lim(C) = C$ implies
that every element in $B(C)$ is non-negative on $\lim(C)$, i.e.\ 
$\lim(C) \subeq B(C)^\star$. 

Using (\ref{eq:reco}), we see that 
for $x \in B(C)^\star$, $c \in C$ and $f \in B(C)$, we have 
$f(x + c) \geq f(c) \geq \inf f(C),$
so that $x + C \subeq C$. This proves $B(C)^\star \subeq \lim(C)$ 
and hence equality. 

(ii) From the Hahn--Banach--Separation Theorem, we further derive 
\break $\oline{B(C)} = (B(C)^\star)^\star = \lim(C)^\star$. 
If $B(C)$ has non-empty interior, it coincides with the interior of its 
closure (\cite{Bou07}, Cor.~II.2.6.1). 

(iii) If $C$ has an equicontinuous weak-$*$-compact 
base $K$ and $x \in V$ satisfies
$\la K, x \ra > \eps$, then $\eps \hat K$ (where $\hat K$ denotes the 
polar of $K$) is a $0$-neighborhood 
in $V$ with $x + \eps \hat K \subeq C^\star$, showing that 
$C^\star$ has interior points. 

If, conversely, $C^\star$ has an interior point $x_0$ and 
$U$ is a convex symmetric $0$-neighborhood with $x_0 + U \subeq C^\star$, 
then the polar set 
$\hat U$ is a weak-$*$-compact equicontinuous subset containing 
$K := \{ \alpha \in C \: \alpha(x_0) = 1\}$. Therefore 
$K$ is  weak-$*$-compact and equicontinuous with $C = \R_+ K$ and 
$0 \not\in K$, i.e., 
$K$ is a base of $C$.
\end{proof}

If $B(C)$ has interior points, the following proposition shows that 
we can reconstruct $C$ from the values of $s_C$ on the open set $B(C)^0$. 

\begin{proposition} \label{prop:3.4} 
Let $V$ be a locally convex space and $C \subeq V'$ be a weak-$*$-closed 
convex subset for which the cone $B(C)$ has interior points. Then 
$$ C = \{ \alpha \in V' \: (\forall x \in B(C)^0)\ \alpha(x) \geq 
\inf \la C, x\ra\}. $$
\end{proposition}

\begin{proof} 
Let 
$D := \{ \alpha \in V' \: (\forall x \in B(C)^0)\ \alpha(x) \geq 
\inf \la C, x\ra\}.$
Then we have $C \subeq D$ and both are weak-$*$-closed convex sets. 
If $\alpha \in D \setminus C$, then (\ref{eq:reco}) implies the existence 
of some $x \in B(C)$ with $\alpha(x)  < \inf \la C, x \ra$. 

Let $x_0 \in B(C)^0$. Then, for each $t \in ]0,1[$, the 
element $x_t := (1-t) x_0 + t x$ is contained in 
$B(C)^0$, so that 
$\alpha(x_t) \geq \inf \la C, x_t \ra.$
Now $F \: [0,1] \to \R, F(t) := -\inf \la C, x_t \ra$ is a lower semicontinuous 
convex function on a real interval, hence continuous (\cite{Ne99}, Cor.~V.3.3). 
Therefore 
$$ \inf \la C, x\ra = \lim_{t \to 1} \inf \la C, x_t \ra 
\leq \lim_{t \to 1} \alpha(x_t) = \alpha(x), $$
and we get $\alpha(x) \geq \inf \la C, x \ra$. 
This contradiction implies that $C = D$. 
\end{proof}

The following lemma is obvious: 
\begin{lemma} \label{lem:3.8} For 
a weak-$*$-closed convex subset $C \subeq V'$ the 
following are equivalent
\begin{description}
\item[(1)] $C$ is weak-$*$-bounded.  
\item[(2)] $B(C) = V$. 
\item[(3)] The polar set $\hat C := \{ v \in V \: |\la C,v\ra| \leq 1\}$ 
is absorbing. 
\end{description}
\end{lemma}

\subsection*{Convex functions on domains in locally convex spaces} 

\begin{definition} Let $X$ be a topological space. 
We say that a real-valued function $f \: X \to \R_\infty := \R\cup\{\infty\}$ 
is {\it lower semicontinuous} if for each 
$x_0 \in X$ and $c < f(x_0)$ there exists a neighborhood $U$ of $x_0$ 
with $\inf f(U) > c$. 
We call $f$ {\it upper semicontinuous} if for each 
$x_0 \in X$ and $d > f(x_0)$ there exists a neighborhood $U$ of $x_0$ 
with $\sup f(U) < d$. 
\end{definition}

\begin{remark} A function $f \: X \to \R_\infty$ 
is lower semicontinuous if and only if its epigraph  
$\epi(f) := \{ (x,t) \in X \times \R \: f(x) \leq t \}$ 
is a closed subset of $X \times \R$. 
\end{remark}

\begin{proposition} \label{prop:3.9} 
{\rm (cf.\ \cite{Bou07}, Prop.~II.2.21)} 
Let $\Omega \subeq V$ be an open convex subset 
and $f \: \Omega \to \R$ a lower semicontinuous convex function. 
Then the following are equivalent
\begin{description}
\item[(1)] $f$ is continuous. 
\item[(2)] $f$ is locally bounded. 
\item[(3)] $f$ is bounded in a neighborhood of one point. 
\item[(4)] $f$ is upper semicontinuous. 
\end{description}
\end{proposition}

\begin{proof} Since $f$ is assumed to be lower semicontinuous, (1) and (4) 
are equivalent. Clearly, (1) $\Rarrow$ (2) $\Rarrow$ (3). 

(3) $\Rarrow$ (2): Suppose 
that $f \leq M$ holds on the open convex neighborhood $U$ of 
$c_0 \in \Omega$. Let $c \in \Omega$. Then there exists an element 
$c_1 \in \Omega$ and $0 < t < 1$ with 
$c = (1 - t)c_1 + t c_0$. Then 
$(1-t) c_1 + t U$ is an open subset of $\Omega$ containing $c$, 
and on this subset we have 
$$ \sup f\big((1-t) c_1 + t U\big) 
\leq (1-t) f(c_1) + t \sup f(U) < \infty. $$
Therefore $f$ is locally bounded. 

(2) $\Rarrow$ (4): Let $c \in \Omega$ and $U$ a closed convex $0$-neighborhood 
in $V$ with $c + U \subeq \Omega$ on which $f$ is bounded and $d > f(c)$. 
Then for each $t \in [0,1]$, we have 
$$ \sup f(c + tU) =  \sup f((1-t)c + t(c + U)) 
\leq (1-t) f(c) + t\sup f(c+U), $$
so that for some $t > 0$ close to $0$, we have $\sup f(c + tU) \leq d$. 
Therefore $f$ is upper semicontinuous in $c$. 
\end{proof}

\begin{remark} Under the circumstances of Lemma~\ref{lem:3.8}, 
$s_C(x) := -\inf \la C, x \ra$
defines a convex function on all of $V$. 
In view of the preceding proposition, this function is locally 
bounded if and only if it is continuous, and this is equivalent 
to the polar set $\hat C$ being a $0$-neighborhood. 
The polar set $\hat C$ is a {\it barrel}, i.e., a closed absolutely convex 
absorbing set. According to the Bipolar Theorem, each barrel 
$B$ coincides with the polar $\hat C$ of its polar $C := \hat B$. 

A locally convex space $V$ is said to be {\it barrelled} if all barrels 
in $V$ are $0$-neighborhoods. In view of the preceding remarks, this means that 
all functions $s_C$ are continuous. The functions $s_C$, $\hat C$ a barrel, 
are precisely the lower semicontinuous  
seminorms on $V$, so that $V$ is barrelled if and  
only if all lower semicontinuous seminorms are continuous 
(cf.\ \cite{Bou07}, \S III.4.1).
\end{remark} 

The following theorem extends this remark to general convex functions. 

\begin{theorem} \label{thm:autcont} Let 
$V$ be a barrelled space, $\Omega \subeq V$ an open 
convex set and $f \: \Omega \to \R$ a lower semicontinuous function. 
Then $f$ is continuous. 
\end{theorem}

\begin{proof} Pick $x_0 \in \Omega$ and let 
$U$ be a closed absolutely convex $0$-neighborhood with 
$x_0 + U \subeq \Omega$. We consider the set 
$$ B := \{ v \in U \: f(x_0 \pm v) \leq f(x_0) + 1\} $$
and claim that $B$ is a barrel. 
Since $f$ is lower semicontinuous, $B$ is a closed convex subset of $U$. Moreover, 
$v \in B$ and $\lambda \in \R$ with $|\lambda| \leq 1$ implies 
$\lambda v \in B$ because 
$$ f(x_0 \pm \lambda v) \leq \conv\{ f(x_0\pm v)\} \leq f(x_0) + 1. $$
We conclude that $B$ is absolutely convex. To see that $B$ is absorbing, 
let $v \in U$ and observe that the lower semicontinuous function 
$h(t) := f(x_0 + tv)$ on $[-1,1]$ 
is continuous (\cite{Ne99}, Cor.~V.3.3). 
Hence there exists a $\mu > 0$ with $\mu v \in B$. This proves 
that $B$ is a closed absolutely convex absorbing set, hence a barrel. 
Since $V$ is barrelled, $B$ is a $0$-neighborhood, and thus 
$f$ is bounded on a neighborhood of $x_0$. The continuity of $f$ now 
follows from Proposition~\ref{prop:3.9}. 
\end{proof}

\begin{remark} For a locally convex space we have the following implications 
$$ \mbox{Banach} \quad \Rarrow \quad \mbox{Fr\'echet} \quad \Rarrow \quad \mbox{Baire} \quad \Rarrow 
\quad \mbox{barrelled}, $$
so that in particular all Fre\'chet spaces are barrelled. 

Furthermore, each locally convex direct limit of barrelled spaces is barrelled, 
which implies that there are barrelled spaces which are not Baire, f.i. 
$V := \R^{(\N)}$, endowed with the finest locally convex topology, 
is such a space. 
\end{remark}

\begin{example} (A barrel which is not a $0$-neighborhood) 
Let $X = c_0$, which is a non-reflexive Banach space and 
$V := \ell^1$ its topological dual, endowed with the 
weak-$*$-topology. Then the closed unit ball $B \subeq V$ 
is a barrel which is not a zero neighborhood because each $0$-neighborhood  
contains a subspace of finite codimension. 
%
\end{example}

\begin{proposition} \label{prop:wlc} Let $C \subeq V'$ be a 
non-empty weak-$*$-closed convex subset and $x \in B(C)$ such that 
the support function $s_C$ is bounded on some neighborhood of $x$. 
Then the following assertions hold: 
\begin{description}
\item[\rm(1)] For each $m \in \R$ the subset 
$C_m := \{ \alpha \in C  \: \alpha(x) \leq m \}$ 
is equicontinuous and weak-$*$-compact. 
\item[\rm(2)] The function $\eta(x) \: C \to \R, 
\eta(x)(\alpha) := \alpha(x)$ is proper. 
\item[\rm(3)] There exists an extreme point $\alpha \in C$ 
with $\alpha(x) = \min \la C, x\ra$. 
\item[\rm(4)] $C$ is weak-$*$-locally compact. 
\end{description}
\end{proposition}

\begin{proof} (1) Pick a $0$-neighborhood $U \subeq V$ for 
which $s_C$ is bounded on $x + U$ by some constant $M$. 
Then $\la C, x + U\ra \geq - M$, and hence 
$$ \la C_m, U \ra \geq \la C_m, x + U \ra - m \geq - M - m. $$
This implies that $s_{C_m}$ is bounded from below on $U$, and 
hence $s_{C_m}$ is bounded from above on $-U$. Therefore 
the polar $\hat C_{m}$ contains a multiple of $U \cap -U$, hence 
is a neighborhood of $0$. This is equivalent to $C_m$ being equicontinuous.
Now the Banach--Alaoglu--Bourbaki Theorem implies that 
$C_m$ is weak-$*$-compact because it is a closed subset of the 
polar set of a $0$-neighborhood in~$U$. 

(2) follows immediately from (1). 

(3) Pick $M > \inf \la C, x\ra$. Then the weak-$*$-compactness of $C_M$ 
implies the existence of a minimal value $m = \min \eta(x)(C)$. 
Then $C_m := \eta(x)^{-1}(m) \cap C$ is a weak-$*$-compact convex set,
so that the Krein--Milman Theorem implies the existence of an extreme point 
$e$ of $C_m$. Since $C_m$ is a face of $C$, $e$ also is an extreme point of $C$.  

(4) For any $\alpha \in C$, (1) implies that the set $C_{\alpha(x)+1}$ is a compact 
neighborhood of $\alpha$ in $C$. Hence $C$ is weak-$*$-locally compact. 
\end{proof}

\begin{remark} If $C \subeq V'$ is a weak-$*$-closed convex subset which is 
locally compact with respect to the weak-$*$-topology, 
then its recession cone
$\lim(C)$ is also locally compact because, for 
each $\alpha \in C$, the subset $\alpha + \lim(C)$ of $C$ is closed. 
Therefore Exercise II.7.21(a) in \cite{Bou07} implies that 
the cone $\lim(C)$ has a weak-$*$-compact base $K$, but if $K$ is not equicontinuous, 
this does not imply that the dual cone $\lim(C)^\star = \oline{B(C)}$ has interior 
points (Lemma~\ref{lem:bcdual}). 
\end{remark}

\begin{remark} If $C \subeq V'$ is weak-$*$-locally compact, then we consider 
the commutative $C^*$-algebra $\cA := C_0(C)$. 
Clearly, the map 
$$ \eta \: V \to M(C_0(C)) \cong C_b(C), \quad 
\eta(x)(\alpha) := e^{i\alpha(x)} $$
is a multiplier action of the abelian topological group 
$V$ on $C_0(C)$. This action is strictly continuous 
if and only if $\eta$ is continuous with respect to the 
compact open topology on $C_b(C)$ (cf.\ \cite{Br77}, Lemma~3.5). 
If this is the case, then the additive map 
$\tilde\eta \: V \to C_b(C), \tilde\eta(x)(\alpha) = \alpha(x)$ 
is also continuous with 
respect to the compact open topology, and this is equivalent 
to the equicontinuity of each compact subset of $C$. 
Then $\eta \: V \to C_b(C)$ defines a strictly continuous multiplier 
action and Example~\ref{exs:hosts}(b) implies that 
$(C_0(C),\eta)$ is a host algebra for $V$. 
\end{remark}

\begin{problem} (a) Suppose that for some $x \in V$ all sets 
$$C_m := \{ \alpha \in C \: \alpha(x) \leq m\}$$ 
are equicontinuous, hence 
in particular weak-$*$-compact. Then $\eta(x)(\alpha) := \alpha(x)$ defines 
a proper function on $C$, showing that $C$ is locally compact. 
The equicontinuity of the sets $C_m$ implies that the cone $\lim(C)$ has 
an equicontinuous weak-$*$-compact base, so that $\lim(C)^\star$ has 
interior points (Lemma~\ref{lem:bcdual}(iii)). Does this imply 
that $B(C)$ has interior points and that 
the support function $s_C$ is bounded on some neighborhood of $x$? 

(b) If, in addition, some function 
$e^{-\tilde\eta(x)}$, $\tilde\eta(x)(\alpha) = \alpha(x)$, 
is contained in $C_0(C)$, then 
$\tilde\eta(x)$ is proper and bounded from below. 
Does the requirement that 
$$ \eta^{-1}(C_0(C)) \subeq B(C) $$
has interior points imply that $s_C$ is bounded on some open set? 
According to Theorem~\ref{thm:autcont}, this is the case if $V$ is barrelled. 
\end{problem}

\section{Host $C^*$-algebras coming from tubes} 

In this section we briefly take a closer look at the host algebras 
of a locally convex space $V$, defined by a complex involutive 
semigroup of the type $S = V + i W$ for an open convex cone 
$W \subeq V$ (Proposition~\ref{prop:hostexam}). In view of 
Remark~\ref{rem:adjoint}(b) and the discussion in Example~\ref{ex:cstarsem}, 
any locally bounded absolute value $\alpha$ 
on such a semigroup leads to the same $C^*$-algebra as an absolute value 
of the form 
$$ \alpha_C(x + iy) := e^{-\inf\la C, y \ra}, $$
where $C \subeq V'$ is a weak-$*$-closed convex subset with 
$W \subeq B(C)$ whose support function 
is locally bounded on $W$. The following theorem provides a converse:

\begin{theorem} Let $C \subeq V'$ be a weak-$*$-closed convex subset 
for which the cone $B(C)$ has interior point and the support function
$s_C$ is locally bounded on $B(C)^0$. Then the following assertions 
hold: 
\begin{description}
\item[\rm(1)] $S := V + i B(C)^0$ is 
an open complex involutive subsemigroup of $V_\C$ and 
$$ \alpha_C(x + iy) := e^{- \inf \la C, y\ra} $$
is a locally bounded absolute value on $S$. 
\item[\rm(2)] $C$ is weak-$*$-locally compact. 
\item[\rm(3)] The map 
$$ \gamma \: S \to C_0(C), \quad 
\gamma(s)(f) := e^{if(s)} $$
induces an isomorphism of $C^*$-algebras 
$C^*(S,\alpha_C) \to C_0(C).$
\item[\rm(4)] The homomorphism $\eta \: V\to C(C,\T) = U(M(C_0(C))), 
\eta(x)(f) := e^{if(x)}$ defines a host algebra of $V$ with a 
strictly smooth multiplier action. The corresponding momentum set is 
$I_\eta = C$. 
\end{description}
\end{theorem}

\begin{proof} (1) is an immediate consequence of the definition. 

 (2) This follows from Proposition~\ref{prop:wlc}. 

(3) For each element 
$s = x + iy \in S$, the continuous function 
$\gamma(s)$ on $C$ satisfies 
$$ \|\gamma(s)\| 
= \sup_{f \in C} |e^{if(s)}| 
= \sup_{f \in C} e^{-f(y)} 
= e^{-\inf \la C, y \ra} = \alpha_C(s). $$

Moreover, for each $\eps > 0$, the subset 
$$ \{ f \in C \: e^{-f(y)} = |\gamma(s)(f)| \geq \eps \} 
= \{ f \in C \: f(y) \leq -\log \eps \} $$
is weak-$*$-compact (Proposition~\ref{prop:wlc}), so that 
$\gamma(s) \in C_0(C)$. 

We thus obtain a  morphism $\gamma \: S \to C_0(C)$ of involutive semigroups 
with $\|\gamma(s)\| \leq \alpha_C(s)$. Hence $\gamma$ is locally bounded, 
and to see that it is holomorphic, it suffices to verify its holomorphy 
on the intersection of $S$ with each complex subspace $E_\C \subeq V_\C$, 
where $E \subeq V$ is finite dimensional. 

Using that the $C^*$-algebra $C_0(C)$ has a realization as a closed 
subalgebra of some algebra of the form $B({\cal H})$, we first 
use \cite{Ne99}, Thm.~VI.2.3, to see that 
$\rho:= \gamma\res_{i B(C)^0 \cap E_\C}$ 
is a norm continuous morphisms of semigroups. Further, 
Proposition~VI.3.2 in \cite{Ne99} implies that  $\rho$ 
extends to a unique holomorphic homomorphism 
$\hat\rho \: S \cap E_\C \to B({\cal H})$, but since $iB(C)^0 \cap E$ 
is totally real in $S \cap E_\C$, the values of the 
unique holomorphic extension $\hat\rho$ also lie in the closed 
subspace $C_0(C)$ of $B({\cal H})$. For each $f \in C$, the function 
$\hat\rho(s)(f)$ on $S \cap E_\C$ is the unique holomorphic extension 
of the function $\eta(s)(f) = \rho(s)(f)$, and from the 
holomorphy of 
$S \to \C, s \mapsto \gamma(s)(f)$ 
we derive that 
$\hat\eta = \gamma\res_{S \cap E_\C}$. This proves that $\gamma$ is holomorphic 
on $S \cap E_\C$, and hence that $\gamma$ is holomorphic. 

Now the universal property of the $C^*$-algebra $C^*(S,\alpha_C)$ leads 
to a unique morphism $\hat\gamma \: C^*(S,\alpha_C) \to C_0(C)$ 
of $C^*$-algebras with 
$\hat\gamma \circ \eta = \gamma$ (Theorem~\ref{thm:1.7}). 

To see that $\hat\gamma$ is injective, we recall that the characters of 
$C^*(S,\alpha_C)$ separate the points, so that it suffices to show that they 
are all of the form $f \mapsto \hat\gamma(f)(f)$ for some $f \in C$. 
Any non-zero character $\chi$ of $C^*(S,\alpha_C)$ is uniquely determined 
by the holomorphic character $\chi_S := \chi \circ \eta \: S \to \C$. 
We claim that $\chi_S(S) \subeq \C^\times$. 
Indeed, $\chi_S(s)=0$ implies $\chi_S(s + S) = \{0\}$, so that 
$\chi_S = 0$ would follow by analytic continuation. 
Now $\chi_S(S) \subeq \C^\times$ shows that we have a corresponding 
smooth character $\chi_V \: V \to \T$, obtained from the smooth 
multiplier action of $V$ on $S$. We further derive 
that there exists some $\beta \in V'$ with
$$ \chi_S(s) = e^{i\beta(s)} \quad \mbox{ for } \quad s \in S= V + i B(C)^0. $$
Since any morphism of $C^*$-algebras is contractive, 
we get 
$$ |\chi_S(x+iy)| = e^{-\beta(y)} \leq \alpha_C(s) = e^{-\inf \la C,y\ra}, $$
i.e., 
$\beta(y) \geq \inf \la C, y \ra$ for $x \in B(C)^0.$
Now we apply Proposition~\ref{prop:3.4} to obtain 
$\beta \in C$. If, conversely, $\beta \in C$, then 
$e^{i\beta}$ defines an $\alpha_C$-bounded holomorphic character of 
$S$, and the universal property of $C^*(S,\alpha_C)$ implies that this 
character extends to a character of $C^*(S,\alpha_C)$. These arguments show 
that the characters of $C^*(S,\alpha_C)$, resp., the 
$\alpha_C$-bounded holomorphic characters of $S$, are of the form 
$s \mapsto \gamma(s)(\beta)$ for some $\beta\in C$. 
As we have already observed above, this implies that $\hat\gamma$ 
is injective, hence an isometric embedding
(\cite{Dix64}, Cor.\ I.8.3). 

In view of the Stone--Weierstra\ss{} Theorem, the fact that 
the functions in $\gamma(S)$ have no zeros and separate the points of 
$C$ implies that $\hat \gamma$ has dense image. We know already that 
$\hat\gamma$ is isometric, so that its range is closed. Therefore 
$\hat\gamma$ is an isomorphism of $C^*$-algebras. 

 (4) First we combine 
Proposition~\ref{prop:holhostalg} with Example~\ref{ex:tube} to 
see that $(C_0(C),\eta)$ is a host algebra of $V$ with strictly 
smooth multiplier action. 

To calculate the corresponding momentum set, we first recall from (3) that 
$S_\alpha = C$, so that the character $\chi_f(\xi) := \xi(f)$ 
defined by $f \in C$ defines a smooth state of $C_0(C)$ with 
$\tilde\chi_f(\eta(x)) = e^{if(x)}$, which leads to 
$\Psi_\eta(\chi_f) = f$, and thus $C \subeq I_\eta$. 
On the other hand, Proposition~\ref{prop:holext} shows that for each 
$y \in B(C)^0$, we have 
$$ e^{-\inf I_\eta(y)} = \|\hat\eta_y(i)\| = \alpha(iy) 
= e^{-\inf \la C, y\ra}, $$
so that $\inf \la C, y \ra = \inf I_\eta(y)$, which leads to 
$I_\eta \subeq C$ (Proposition~\ref{prop:3.4}). 
\end{proof} 

The preceding theorem implies in particular that each weak-$*$-closed convex 
subset $C \subeq V'$ whose support function is bounded on some open 
subset actually occurs as the momentum set of some host algebra of $V$. 
Conversely, we have seen in Example~\ref{ex:cstarsem}(c) that 
all host algebras defined by complex semigroups of the form 
$V + i W$, $W$ an open convex cone in $V$, are of this form. We thus obtain a 
complete picture for the case where $G = (V,+)$ is the additive group 
of a  locally convex space. 

\begin{remark} One can also develop a holomorphic representation 
theory of tubes of the form $V + i W$ by starting with representations 
of the cone $iW \cong W$ by selfadjoint operators. This program 
has been carried out in great generality by H.~Gl\"ockner in 
\cite{Gl03} (cf.\ also \cite{Gl00}). 
\end{remark}

\section{The finite dimensional case} 

In the preceding section we have described all host algebras of 
locally convex spaces defined by complex host semigroups. 
In the non-commutative case this turns out to be much harder. However, 
using \cite{Ne99}, we also obtain a complete picture for 
finite dimensional groups. We shall take up the investigation of the 
non-abelian infinite dimensional case in the future. 

We write  $S^\infty(G)$ for the set of {\it smooth states of $G$}, i.e., 
the set of smooth positive definite functions normalized by 
$\phi(\1) = 1$. We consider $S^\infty(G)$ as a convex subset of the set 
$S(G)$ of continuous states of $G$, which in turn can be identified 
with the state space $S(C^*(G))$ of the group algebra $C^*(G)$. 
We recall the following result from \cite{Ne99}, Prop.~X.6.17. 

\begin{proposition}   \label{prop:10.6.17} 
Let $C \subeq \L(G)^*$ be a closed convex invariant subset 
and 
$$ \ev \: S^\infty(G) \to \L(G)^*, \quad \phi \mapsto 
\frac{1}{i} d\phi(\1). $$
Then the annihilator $I_C := \ev^{-1}(C)^\bot$
is an ideal of $C^*(G)$. 
The non-degenerate representations of the quotient algebra 
$C^*(G)_C := C^*(G)/I_C$ correspond to those continuous unitary
representations $(\pi, {\cal H})$ of $G$ satisfying $I_\pi \subeq C$. 
\end{proposition}

\begin{theorem} Let $G$ be a connected finite dimensional Lie group, 
$(S,\eta_S, W)$ a host semigroup of $G$ and $\alpha$ a locally bounded 
absolute value on $S$. Then the following assertions hold: 
\begin{description}
\item[\rm(a)] The host algebra $(C^*(S,\alpha),\eta)$ 
is a quotient of $C^*(G)$. 
\item[\rm(b)] If, in addition, $G$ acts on $C^*(S,\alpha)$ with discrete 
kernel and the polar map 
$G \times W \to S, (g,x) \mapsto g\Exp x$ is a diffeomorphism, 
then $C^*(S,\alpha) \cong C^*(G)_{I_\eta}$. 
\end{description}
\end{theorem}

\begin{proof} (a) Let $(\pi, {\cal H})$ be the universal 
representation of $C^*(S,\alpha)$. Then we have a holomorphic 
representation $\hat\pi \: S \to B({\cal H})$ whose image generates 
the $C^*$-algebra $\cB := \pi(C^*(S,\alpha)) \cong C^*(S,\alpha)$ 
(Theorem~\ref{thm:1.7}). 

Let $\pi_G \: G \to U({\cal H})$ denote the corresponding 
unitary representation of $G$ and $\hat\pi_G$ the associated 
representation of $C^*(G)$. Since $G$ acts smoothly by multipliers on $S$, 
we obtain a continuous multiplier action of $G$ on $\cB$, and this leads to 
$\hat\pi_G(C^*(G)) \cB \subeq \cB,$
where the left hand side is dense in $\cB$. 

On the other hand, $G$ also acts continuously by unitary multipliers on 
$C^*(G)$, hence on $\hat\pi_G(C^*(G))$. For each $x \in W$, 
we now obtain a morphism of 
$C^*$-algebras $\tilde\pi_x \: 
C^*(\R) \to M(\hat\pi_G(C^*(G)))$ which factors through 
some quotient $C_0([m,\infty[)$. This implies that 
$\hat\pi(\Exp(x)) \in M(\hat\pi_G(C^*(G)))$, and from that 
we obtain $\hat\pi(S) \subeq M(\hat\pi_G(C^*(G)))$ 
by analytic continuation and (HS3) in the definition of a host semigroup. 
This in turn leads to $\cB\hat\pi_G(C^*(G)) \subeq \hat\pi_G(C^*(G))$.  
Since $\tilde\pi_x(C_0([m,\infty[))\hat\pi_G(C^*(G))$ is dense in 
$\hat\pi_G(C^*(G))$, we see that 
$\cB\hat\pi_G(C^*(G))$ spans a dense subspace of $\hat\pi_G(C^*(G))$. 
We now arrive that 
$$ \cB = \oline{\Spann(\hat\pi_G(C^*(G))\cB)}= \hat\pi_G(C^*(G)),  $$
and this proves (a). 

(b) From the proof of 
(1) $\Leftrightarrow$ (2) in Proposition~\ref{prop:holext} 
we recall that $I_\pi = I_\eta$. Let $C := I_\eta$. Then the ideal 
$I_C$ of $C^*(G)$ annihilates all states of the form $\pi^v$, 
$v \in {\cal H}$, so that $I_C \subeq \ker \hat\pi_G$. 

By assumption, $d\pi$ is a faithful representation of $\L(G)$, so that 
$\{0\} = \ker d\pi = I_\pi^\bot$ implies that $I_\pi$ spans the dual space 
$\L(G)'$. Therefore the open cone $W \subeq B(I_\pi)$ satisfies 
$\Spec(\ad x) \subeq i \R$ for each $x \in W$ 
(\cite{Ne99}, Prop.~VII.3.4(b)), i.e., $W$ is a weakly elliptic cone 
(\cite{Ne99}, Def.~XI.1.11). 
The construction in Section~XI.1 in \cite{Ne99} now leads to a 
complex involutive semigroup $\Gamma_G(W)$ for which the polar 
map 
$$ G \times W \to \Gamma_G(W), \quad (g,x) \mapsto g \Exp(x) $$
is a diffeomorphism. 
According to the Holomorphic Extension Theorem 
(\cite{Ne99}, XI.2.3), each unitary representation 
$(\rho, {\cal K})$ of $G$ with $I_\rho \subeq C = I_\eta$ extends via 
$$ \rho \: \Gamma_G(W) \to \hat\pi(S) \subeq \cB \subeq B({\cal H}), \quad 
g \Exp x \mapsto \pi(g) e^{i\cdot d\pi(x)} $$
to a holomorphic representation $\hat\rho$ of $\Gamma_G(W)$ with 
$$\|\hat\rho(g \Exp x)\| 
= e^{-\inf \la I_\rho, x\ra} 
\leq e^{-\inf \la C, x\ra} = e^{-\inf \la I_\eta, x\ra} 
\leq \alpha(g \Exp x), $$
hence to a representation of $\cB$. In view of (a), these representations 
separate the points of $C^*(G)_C$, which implies that 
$I_C = \ker \hat\pi_G$. We finally obtain 
$C^*(S,\alpha) \cong \cB \cong C^*(G)_C$. 
\end{proof}

\begin{remark} \label{rem:9.6.25}
If $G = V$ is a finite dimensional vector space, $W \subeq V$ an open convex 
cone, and 
$C \subeq V'$ a closed convex subset, then 
$C^*(G) \cong C_0(\hat G) \cong C_0(V')$, and 
the definition of $C^*(G)_C$, implies that 
$C^*(G)_C \cong C_0(C)$. We have already seen in Section~5 that 
this is $C^*(S, \alpha)$ for $S = V + i W$ and 
$\alpha(x + iy) = e^{-\inf \la C, y \ra}$. 
\end{remark}

\section{Open Problems}

\begin{problem} (Invariant convex geometry of Lie algebras) 
Study open invariant convex cones in the Lie algebra 
$\L(G)$ of an infinite dimensional Lie group $G$. Here are some concrete 
problems: 
\begin{itemize}
\item Does it have any consequence for the spectrum of $\ad x$ if $x$ 
is contained in an open invariant convex cone $W$ not containing 
affine lines? What can be said about the stabilizer of $x$ in $G$ and its 
action on the Lie algebra $\L(G)$? In the finite-dimensional 
case it acts like a compact group and, consequently, $\ad x$ is semisimple 
with $\Spec(\ad x) \subeq i \R$ (cf.\ \cite{Ne99}). 

\item Develop a structure theory for coadjoint orbits 
${\cal O}_f := \Ad^*(G).f \subeq \L(G)'$ for which the 
weak-$*$-closed convex hull $C_f$ has the property that 
$B(C_f)$ has interior points and the support function $s_{C_f}$ is 
locally bounded on the interior. For any such orbit which separates the 
points of $\L(G)$ (which can always be arranged after factorization of a 
closed ideal), the open cone $B(C_f)^0$ does not contain 
affine lines. It is a natural question under which circumstances the 
coadjoint orbit is closed. The geometric setup leads to the alternative 
that either ${\cal O}_f$ consists of extreme points of its 
weak-$*$-closed convex hull or not, where the latter case does not 
arise for closed orbits in the finite-dimensional case 
(cf.\ Section VIII.1 in \cite{Ne99}). 
\end{itemize}
\end{problem}

\begin{remark} If $x \in B(C)^0$ and $f \in C$ is a unique minimum of the 
function $\tilde\eta(x)(\alpha) = \alpha(x)$ on $C$, then the stabilizer 
of $x$ in $G$ is contained in the stabilizer of $G_f$ and it also preserves 
all weak-$*$-compact subsets $$C_m = \{ \alpha \in C \: \alpha(x) \leq m\}$$ 
of $\L(G)'$. This situation should lead to interesting geometric structures 
on the coadjoint orbit ${\cal O}_f$, such as weak K\"ahler structures. 
\end{remark}

\begin{example} Let $\cA$ be a unital $C^*$-algebra and $G = U(\cA)$ its 
unitary group, considered as a Banach--Lie group. Then for each state 
$\phi \in S(\cA)$ the functional 
$$-i\phi \: \L(G) = \{ a \in \cA \: a^* = - a\} = \fu(\cA) \to \R$$ 
is real-valued. If $\phi$ is a pure state, i.e., an extreme point of 
$S(\cA)$, then the coadjoint orbit ${\cal O}_{-i\phi}$ consists of 
extreme points of its weak-$*$-closed convex hull, which is the 
weak-$*$-closed convex face of $-i S(\cA)$, generated by $-i\phi$ 
(\cite{Ne02}, Thm.~III.1). 
\end{example}

\begin{problem} (Holomorphic extensions of unitary representations) 
  \begin{enumerate}
  \item[\rm(1)] Suppose that $(\pi, {\cal H})$ 
is a unitary representation of the infinite dimensional Lie group 
$G$ for which the subspace ${\cal H}^\infty$ of smooth vectors is 
dense and $B(I_\pi)$ has interior points. Let $x \in B(I_\pi)^0$. 
Is the selfadjoint operator $e^{i\cdot d\pi(x)} \in B({\cal H})$ a 
smooth vectors for the multiplier action of $G$ on $B({\cal H})$? 
   \item[\rm(2)] Is $\pi(G) e^{i\cdot d\pi(B(I_\pi)^0)}$ 
a subsemigroup of $B({\cal H})$? 
   \item[\rm(3)] Suppose that there exists a host semigroup $(S,\eta,W)$ for $G$ 
for which the polar map $G \times W \to S, (g,x) \mapsto g \Exp(x)$ 
is a diffeomorphism. Assume that $W \subeq B(I_\pi)$ for some smooth unitary 
representation $(\pi, {\cal H})$ of $G$. Is the map 
$\hat\pi \: S \to B({\cal H}), g\Exp(x) \mapsto \pi(g)e^{i\cdot d\pi(x)}$ 
a holomorphic representation? 
   \item[\rm(4)] If $\L(G)$ contains a dense locally finite subalgebra, 
many of the arguments seem to be reducible to the finite dimensional situation.
  \end{enumerate}
\end{problem}

\begin{problem} (Existence of complex semigroups) 
Let $G$ be a Banach--Lie group (or locally exponential) and 
$W \subeq \L(G)$ an open convex cone satisfying 
$\Spec(\ad x) \subeq i \R$ for each $x \in W$. 
We assume that $G$ has a faithful universal complexification 
$\eta \: G \to G_\C$ (which is locally exponential if $G$ is not Banach). 
Is it true that $G \exp(iW)$ a subsemigroup of $G_\C$? 
For some interesting examples of such semigroups we refer to \cite{Ne01}. 
\end{problem}

\section{Appendix: Some useful facts on multiplier algebras} 

The following results are used in our discussion of general host algebras 
of topological groups. 
\begin{theorem} \label{thm:7.3}
{\rm(\cite{Pa94}, Th.~5.2.2)} Let ${\cal A}$ be a Banach 
algebra with bounded left approximate identity and $T \: \cA \to B(X)$ 
a continuous representation of $\cA$ on the Banach space $X$. 
Then for each $y \in \oline{\Spann(T(\cA)X)}$ there are elements 
$a \in \cA$ and $x \in X$ with $y = T(a)x$. 
\end{theorem}

\begin{corollary} \label{cor:7.6}
If $\cA$ and $\cB$ are $C^*$-algebras and 
$\pi \: \cA \to M(\cB)$ is a homomorphism for which 
$\pi(\cA)\cB$ is dense in $\cB$, then each element $y \in \cB$ can be written 
as $\pi(a)b$ for some $a \in \cA$ and $b \in \cB$. 
\end{corollary}

\begin{proposition} \label{prop:b.2} 
Let $\cA$ and $\cB$ be $C^*$-algebras. 
For each morphism $\alpha \: \cA \to M(\cB)$ of $C^*$-algebras for which 
$\alpha(\cA)\cB$ is dense in $\cB$,  
there exists a unique morphism of $C^*$-algebras 
$\tilde\alpha \: M(\cA) \to M(\cB)$ extending $\alpha$,  and 
$\tilde\alpha$ is strictly continuous. 
\end{proposition}

\begin{proof} The uniqueness of $\tilde\alpha$ follows from 
the density of $\alpha(\cA)\cB$ in $\cB$ and 
$\tilde\alpha(m)\alpha(a)b = \alpha(ma)b$ for 
$m \in M(\cA)$, $a \in \cA$, $b \in \cB$. 

For the existence, we realize $\cB$ as a closed $*$-subalgebra 
of some $B({\cal H})$ for which the representation on 
${\cal H}$ is non-degenerate. Then 
$$M(\cB) \cong \{ T \in \cB({\cal H}) \: T\cB + \cB T \subeq \cB\} $$
(Examples~\ref{ex:multi}) 
and we interprete $\alpha$ as a representation of $\cA$ on ${\cal H}$. 

We claim that $\alpha$ is non-degenerate. Indeed, if 
$\alpha(\cA)v = \{0\}$, then 
$$\{0\} 
= \la \alpha(\cA)v, \cB{\cal H} \ra 
= \la \cB\alpha(\cA)v, {\cal H} \ra $$
implies $\cB\alpha(\cA)v = \{0\}$, and since 
$\cB\alpha(\cA) = (\alpha(\cA)\cB)^*$ is dense in $\cB$, we obtain $v = 0$. 

As $\alpha$ is non-degenerate, there exists a unique extension 
$\tilde\alpha\: M(\cA) \to B({\cal H})$ with 
$\tilde\alpha(m)\alpha(a) = \alpha(ma)$ for 
$m \in M(\cA)$ and $a \in \cA$. 
Then 
$$ \tilde\alpha(m) \alpha(\cA)\cB \subeq \alpha(\cA)\cB \subeq \cB 
\quad \mbox{ and } \quad 
\cB\alpha(\cA) \tilde\alpha(m) \subeq \cB\alpha(\cA) \subeq \cB, $$
and the density of $\alpha(\cA)\cB$, resp., $\cB\alpha(\cA)$ in $\cB$ implies 
that $\tilde\alpha(M(\cA)) \subeq M(\cB)$. 

Suppose that $c_i \to c$ strictly in $M(\cA)$. 
Then we have for 
$a \in \cA$ and $b \in \cB$ the relation 
$\tilde\alpha(c_i) \alpha(a)b = \alpha(c_ia)b \to \alpha(ca)b,$
and since $\cB = \alpha(\cA)\cB$ (Corollary~\ref{cor:7.6}), 
we get $\tilde\alpha(c_i)b \to \tilde\alpha(c)b$ for each $b \in \cB$, 
showing that $\tilde\alpha(c_i) \to \tilde\alpha(c)$ in the 
strict topology on $M(\cB)$. 
\end{proof}

\begin{proposition} \label{prop:b.1} 
Let $\cA$ be a $C^*$-algebra and 
$M(\cA)$ its multiplier algebra. Then the following properties 
of a representation $(\rho, {\cal H})$ of $M(\cA)$ are equivalent
\begin{enumerate}
\item[\rm(1)] $\rho\res_\cA$ is non-degenerate. 
\item[\rm(2)] $\rho = \tilde\pi$ for some non-degenerate 
representation of $\cA$. 
\item[\rm(3)] The representation 
$\rho$ is continuous with respect to the strict topology on $B$ 
and the strong operator topology on $B({\cal H})$. 
\end{enumerate}
\end{proposition} 

\begin{proof}  (1) $\Rarrow$ (2): Let $\pi := \rho\res_\cA$ and 
assume that this representation is non-degenerate, so that it has an 
extension to a representation $\tilde\pi$ of $M(\cA)$ which is uniquely 
determined by $\tilde\pi(m)\pi(a) = \pi(ma)$ for 
$m \in M(\cA)$ and $a \in \cA$. Since $\rho$ also satisfies 
$\rho(m)\pi(a) = \rho(m)\rho(a)= \rho(ma) = \pi(ma)$, we obtain 
$\rho = \tilde\pi$. 

(2) $\Rarrow$ (3): Let $(b_i)_{i \in I}$ be a net in $B$ converging 
to some $b \in B$ with respect to the strict topology, i.e., 
$b_i a \to ba$ and $ab_i \to ab$
for each $a \in \cA$. For each $v \in {\cal H}$ and $a \in \cA$ we then have
$b_i.(a.v) = (b_i a).v \to (ba).v = b.(a.v)$
and 
$$ b_i^*.(a.v) = (a^* b_i)^*.v \to (a^*b)^*.v = b^*.(a.v). $$
We conclude that 
$b_i.v \to b.v$ and $b_i^*.v \to b^*.v$ hold for all vectors $v \in 
\Spann(\cA{\cal H})$, but Theorem~\ref{thm:7.3} implies that 
${\cal H} = \cA{\cal H}$, proving (3). 

(3) $\Rarrow$ (1): Let $(u_i)_{i \in I}$ be an approximate identity 
in $\cA$. Then $u_i \to \1$ in the strict topology on $M(\cA)$. Hence we get for each 
$v \in {\cal H}$ the relation 
$u_i.v \to v$, showing that $\rho\res_\cA$ is non-degenerate. 
\end{proof}

\end{document}